\documentclass[12pt]{amsart}         
\usepackage{amscd}                   
\usepackage{amssymb}

\newtheorem{thm}{Theorem}[section]
\newtheorem{lem}[thm]{Lemma}

\newtheorem{cor}[thm]{Corollary}

\renewcommand{\thestep}{}
\newtheorem{prop}[thm]{Proposition}

\theoremstyle{definition}

\renewcommand{\thecase}{}
\newtheorem{conj}[thm]{Conjecture}

\newtheorem{defn}[thm]{Definition}
\newtheorem{exmp}[thm]{Example}
 
\newtheorem{rmk}[thm]{Remark}

\theoremstyle{remark}

\makeatletter
\def\alphenumi{
  \def\theenumi{\alph{enumi}}
  \def\p@enumi{\theenumi}
  \def\labelenumi{(\@alph\c@enumi)}}
\makeatother

\makeatletter
\def\thecase{\@arabic\c@case}
\makeatother

\numberwithin{equation}{section}

\makeatletter
\def\thestep{\@arabic\c@step}
\makeatother

\newenvironment{pf*}[1]{\begin{proof}[#1]}{\end{proof}}

\newcommand\barB{{\bar B}}

\newcommand\barM{{\bar{M}}}

\newcommand\barV{{\bar{V}}}
\newcommand\barW{{\bar{W}}}

\newcommand\AAA{\mathbb{A}}

\newcommand\CC{\mathbb{C}}

\newcommand\PP{\mathbb{P}}
\newcommand\QQ{\mathbb{Q}}
\newcommand\RR{\mathbb{R}}

\newcommand\ZZ{\mathbb{Z}}

\newcommand\bga{{\boldsymbol{\gamma}}}

\newcommand\bvarphi{{\boldsymbol{\varphi}}}

\newcommand\bE{{\mathbf{E}}}

\newcommand\bL{{\mathbf{L}}}

\newcommand\bM{{\mathbf{M}}}

\newcommand\bx{{\mathbf{x}}}

\newcommand\csC{{{}^\circ\mathcal{C}}}
\newcommand\ssC{{{}^\circ\mathcal{C}}}

\newcommand\ssG{{{}^\circ\mathcal{G}}}
\newcommand\cM{{{}^\circ{M}}}

\newcommand{\cov}{\nabla}

\newcommand{\rd}{\partial}

\newcommand\half{{\textstyle{\frac{1}{2}}}}
\newcommand\third{{\textstyle{\frac{1}{3}}}}

\newcommand\quarter{{\textstyle{\frac{1}{4}}}}
\newcommand\threehalf{{\textstyle{\frac{3}{2}}}}

\newcommand\fm{{\mathfrak{m}}}

\newcommand\fs{{\mathfrak{s}}}
\newcommand\fS{{\mathfrak{S}}}

\newcommand\al{\alpha}

\newcommand\De{\Delta}

\newcommand\eps{\varepsilon}

\newcommand\La{\Lambda}

\newcommand\Om{\Omega}

\newcommand\gl{{\mathfrak{g}\mathfrak{l}}}
\newcommand\fsl{{\mathfrak{s}\mathfrak{l}}}
\newcommand\so{{\mathfrak{s}\mathfrak{o}}}

\newcommand\su{{\mathfrak{s}\mathfrak{u}}}
\newcommand\fu{{\mathfrak{u}}}

\newcommand\PU{\operatorname{PU}}

\newcommand\SO{\operatorname{SO}}

\newcommand\SU{\operatorname{SU}}
\newcommand\U{\operatorname{U}}

\newcommand\less{\setminus}

\newcommand{\8}{\infty}

\newcommand\ad{{\operatorname{ad}}}

\newcommand\asd{{\operatorname{asd}}}
\newcommand\barasd{{\overline{\operatorname{asd}}}}

\newcommand\Coker{\operatorname{Coker}}

\newcommand\End{\operatorname{End}}

\newcommand\Gl{\operatorname{Gl}}

\newcommand\Hom{\operatorname{Hom}}

\newcommand\Ind{\operatorname{Ind}}
\newcommand\Imag{\operatorname{Im}}

\newcommand\Ker{\operatorname{Ker}}

\newcommand\loc{\operatorname{loc}}

\newcommand\red{\operatorname{red}}
\newcommand\Red{\operatorname{Red}}

\newcommand\sign{\operatorname{sign}}

\newcommand\Stab{\operatorname{Stab}}

\newcommand\sw{{\operatorname{sw}}}
\newcommand\SW{{SW}}
\newcommand\Sym{\operatorname{Sym}}
\newcommand\Tor{\operatorname{Tor}}

\newcommand\even{{\mathrm{even}}}

\newcommand\id{{\mathrm{id}}}
\newcommand\odd{{\mathrm{odd}}}

\newcommand\spinc{\text{$\text{spin}^c$ }}
\newcommand\Spinc{\text{$\text{Spin}^c$}}

\newcommand\sA{{\mathcal{A}}}
\newcommand\sB{{\mathcal{B}}}
\newcommand\sC{{\mathcal{C}}}
\newcommand\sD{{\mathcal{D}}}

\newcommand\sG{{\mathcal{G}}}

\newcommand\sN{{\mathcal{N}}}
\newcommand\sO{{\mathcal{O}}}

\newcommand\sV{{\mathcal{V}}}
\newcommand\sW{{\mathcal{W}}}

\newcommand\tsC{{\tilde\sC}}

\newcommand\tM{{\tilde M}}

\newcommand\vecfm{{\vec{\fm}}}
\newcommand\vectau{{\vec{\tau}}}
\newcommand\vecvartheta{{\vec{\vartheta}}}

\begin{document}
\title[$PU(2)$ Monopoles and Relations between 4-Manifold Invariants]
{PU(2) Monopoles and Relations between Four-Manifold Invariants$\text{}^1$}  
\author[Paul M. N. Feehan]{Paul M. N. Feehan}
\address{Department of Mathematics\\
Harvard University\\
One Oxford Street\\
Cambridge, MA 02138}
\email{feehan@math.harvard.edu}
\author[Thomas G. Leness]{Thomas G. Leness}
\address{Department of Mathematics\\
Michigan State University\\
East Lansing, MI 48824--1027}
\email{leness@math.msu.edu}
\thanks{The first author was supported in part by an NSF Mathematical 
Sciences Postdoctoral Fellowship under grant DMS 9306061}
\footnotetext[1]{Slightly revised version to appear in {\em Topology and
its Applications\/}, Proceedings of the Georgia Topology Conference,
Atlanta, GA, June 1996; \texttt{dg-ga/9709022}.}
\subjclass{}
\keywords{}
\maketitle


\section{Introduction}
The principal objective of our series of articles \cite{FL1,FL2,FL3,FL4}
and beyond, for which we provide a brief survey here, is to prove the
analogue of the Kotschick-Morgan conjecture for $\PU(2)$ monopoles
suggested by Pidstrigach and Tyurin \cite{PTLocal}.  This in turn should
lead to a proof of Witten's conjecture concerning the relation between
Donaldson and Seiberg-Witten invariants and a deeper understanding of the
highly successful role of gauge theory in smooth four-manifold topology. We
describe Witten's conjecture below and outline the program (see
\cite{HPP1,HPP2,LabMar1,LabMar2,OTQuaternion,OTVortex,PTCambridge,PTLocal}),
to prove this conjecture 
using $\PU(2)$ monopoles.  While the basic ideas in this program are by now
well-known, the profound analytical difficulties inherent in attempts to
implement it are perhaps much less well-known and so we feel it is
worthwhile to describe some of these analytical problems here. These
analytical difficulties involve the gluing construction of links of
lower-level moduli spaces of $\U(1)$ monopoles contained in the Uhlenbeck
compactification of the moduli space of $\PU(2)$ monopoles.  The question
of existence of perturbations for the $\PU(2)$ monopole equations, yielding
both useful transversality results and an Uhlenbeck compactification for
the perturbed moduli space, is a fairly substantial one in its own right
\cite{FL1}. We describe these transversality and compactness
results here, along with some of our calculations of Donaldson invariants
in terms of Seiberg-Witten invariants from \cite{FL2} and a brief overview of
issues concerning the gluing theory from \cite{FL3,FL4} and its applications.

First, to explain Witten's conjecture we recall that a closed, smooth
four-manifold $X$ is said to have {\em Seiberg-Witten simple type\/} if the
Seiberg-Witten moduli spaces corresponding to non-zero Seiberg-Witten
invariants are all zero-dimensional. The manifold $X$ has 
{\em Kronheimer-Mrowka simple type\/} provided the Donaldson invariants
corresponding to products $z$ of homology classes in $H_\bullet(X;\ZZ)$ and
a generator $x\in H_0(X;\ZZ)$ are related by $D_X^w(x^2z)=4D_X^w(z)$.
Kronheimer and Mrowka \cite{KMStructure} (see also
\cite{FSStructure}) showed that the Donaldson series
of a four-manifold of Kronheimer-Mrowka simple type with $b^1(X)=0$ and odd
$b^+(X)\ge 3$ is given by 
\begin{equation}
\sD^w = e^{Q/2}\sum_{r=1}^s(-1)^{(w^2+wK_r)/2}a_r e^{K_r},
\label{eq:KMFormula}
\end{equation}
where $w$ is a line bundle over $X$, $Q$ is the intersection form on
$H_2(X;\ZZ)$, the coefficients $a_r$ are non-zero rational numbers, and the
$K_r\in H^2(X;\ZZ)$ are the {\em Kronheimer-Mrowka basic classes\/}. Let
$\Spinc(X)$ be the set of isomorphism classes of \spinc structures on
$X$ and let $e(X)$ and $\sigma(X)$ denote the Euler characteristic and
signature of $X$, respectively.

\begin{conj}[Witten]\cite{Witten}
Suppose $X$ is a closed, oriented four-manifold with $b^1(X)=0$ and odd
$b^+(X)\ge 3$, equipped with a homology orientation and a line bundle $w$.
Then $X$ has Kronheimer-Mrowka simple type if and only if it has
Seiberg-Witten simple type. If $X$ has simple type, then the
Kronheimer-Mrowka basic classes are given by
$$
\{c_1(W^+_{\fs}): \fs \in \Spinc(X) \text{ such that } \SW(\fs) \ne 0\},
$$
where $c_1(\fs) := c_1(W^+_{\fs})$ 
and $W^\pm_{\fs}$ are the \spinc bundles
associated to $\fs$ with some choice of Riemannian metric on $X$;
furthermore, the Donaldson series for $X$ is given by 
\begin{equation}
\sD^w
=2^{2+(7e+11\sigma)/4}e^{Q/2}
\sum_{\fs\in\Spinc(X)}(-1)^{(w^2+wc_1(\fs))/2}\SW(\fs)e^{c_1(\fs)}.
\label{eq:WittenFormula}
\end{equation}
\end{conj}

The conjecture holds for all four-manifolds whose Donaldson and
Seiberg-Witten invariants have been independently computed.  The quantum
field theory argument giving the above relation when $b^+(X)\ge 3$ has
recently been extended by Moore and Witten \cite{MooreWitten} to allow
$b^+(X) \ge 1$, $b^1(X)\ge 0$, and four-manifolds $X$ of possibly
non-simple type. The mathematical approach to this conjecture uses a moduli
space of solutions to the $\PU(2)$ monopole equations --- which generalize
the $\U(1)$ monopole equations of Seiberg and Witten --- to construct a
cobordism between links of the compact moduli spaces of $\U(1)$ monopoles
of Seiberg-Witten type and the Donaldson moduli space of anti-self-dual
connections, which appear as singularities in this larger moduli
space. Moreover, this approach should give a precise relation between the
Donaldson and Seiberg-Witten invariants even for four-manifolds not of
simple type. This is an important point since there are no known examples
of four-manifolds with $b^+>1$ violating either of the simple type
conditions, so we would expect to gain a greater understanding of these
conditions from such a general relation.

The moduli space of $\PU(2)$ monopoles is non-compact and has an Uhlenbeck
compactification similar to that of the moduli space of anti-self-dual
connections. The substantial analytical difficulties are due to the
contributions of moduli spaces of $\U(1)$ monopoles (cobordant to standard
Seiberg-Witten moduli spaces) in the lower Uhlenbeck levels --- the
`reducibles' at the boundary of the Uhlenbeck compactification. Many of
these problems had never been resolved even in the case of Donaldson theory
where they arise, albeit in a rather simpler form, in attempts to prove the
Kotschick-Morgan conjecture for Donaldson invariants. The Kotschick-Morgan
conjecture for Donaldson invariants of four-manifolds $X$ with $b^+(X)=1$
asserts that the invariants computed using metrics lying in different
chambers of the positive cone of $H^2(X;\RR)/\RR^*$ differ by terms
depending only the homotopy type of $X$ \cite{KoM}. The heart of the
problem there lies in describing the links of the reducible connections in
the lower Uhlenbeck levels via gluing and then computing integrals of the
Donaldson cohomology classes over those links. To date, links of this type
in anti-self-dual moduli spaces have been described and their pairings with
cohomology classes computed in only a few relatively simple special cases
\cite{DonApplic,DonConn,DonHCobord,DK,Leness,Y}: the methods used there
fall far short of what is needed to complete the $\PU(2)$ monopole program
to prove the equivalence between Donaldson and Seiberg-Witten invariants.
By assuming the Kotschick-Morgan conjecture, G\"ottsche has computed the
coefficients of the wall-crossing formula in \cite{KoM} in terms of modular
forms by exploiting the presumed homotopy invariance of the coefficients
\cite{Goettsche}. A related approach to the Witten conjecture has been
taken so far by Pidstrigach and Tyurin \cite{PTLocal}: they assume a
$\PU(2)$ monopole analogue of the Kotschick-Morgan conjecture and argue
that it can be used to compute the required integrals of analogues of the
Donaldson cohomology classes over the links of the lower-level moduli
spaces of $\U(1)$ monopoles. For a survey of the work of Okonek and Teleman
on non-abelian monopoles, with applications to algebraic geometry and the
conjectured relations between Donaldson and Seiberg-Witten invariants, we
refer the reader to their article \cite{OTSurvey} and the references
contained therein.

In \S \ref{sec:TransvCompact} we describe the $\PU(2)$ monopole equations,
the holonomy perturbations we use in order to achieve transversality, and the
Uhlenbeck compactification for the 
perturbed moduli space of $\PU(2)$ monopoles. In
\S \ref{sec:Cohomology} we describe the cohomology classes, the links of the
moduli spaces of anti-self-dual connections and Seiberg-Witten monopoles
appearing in the top Uhlenbeck level, their orientations, and the
relation between the Donaldson and Seiberg-Witten invariants when 
the moduli spaces of $\U(1)$ monopoles appear only in the top Uhlenbeck
level. Finally, in \S
\ref{sec:GluingKMConjecture} we describe the Kotschick-Morgan conjecture,
its analogue in the case of $\PU(2)$ monopoles and how this might be used
to prove Witten's conjecture. We also describe the need for gluing, survey
some of the results from \cite{FL3, FL4} and describe a few of the more
prominent difficulties which arise when gluing $\PU(2)$ monopoles.
Detailed proofs of all our results appear elsewhere
\cite{FL1,FL2,FL3,FL4}, so we just sketch the main ideas here.

\subsubsection*{Acknowledgements}
The authors warmly thank Gordana Mati\'c (without whose gentle
encouragement and considerable patience this article might not have been
written), the organizers of the 1996 Georgia Topology Conference, and the
Mathematics Department of the University of Georgia, Athens, for their
hospitality. We also thank Peter Ozsv\'ath for many helpful comments.
Finally, we thank the Mathematics Departments
at Harvard and Michigan State University
and the National Science Foundation for their generous support during the
preparation of this article.


\section{Holonomy perturbations,
transversality, and Uhlenbeck compactness}
\label{sec:TransvCompact}
We consider Hermitian two-plane bundles $E$ over $X$ whose determinant line
bundles $\det E$ are isomorphic to a fixed Hermitian line bundle over $X$
endowed with a fixed $C^\8$, unitary connection.  Choose a Riemannian
metric on $X$ and let 
$\fs_0:=(\rho,W)$ be a \spinc structure on $X$, where
$\rho:T^*X\to\End W$ is the Clifford map, and the Hermitian
four-plane bundle $W=W^+\oplus W^-$ is endowed with a $C^\8$ \spinc
connection. {\em The \spinc structure $(\rho,W)$, the \spinc
connection on $W$, and the Hermitian line bundle together with
its connection are fixed once and for all.\/}

Let $k\ge 2$ be an integer and let $\sA_E$ be the space of $L^2_k$
connections $A$ on the $\U(2)$ bundle $E$ all inducing the fixed
determinant connection on $\det E$.  Equivalently, following \cite[\S
2(i)]{KMStructure}, we may view $\sA_E$ be the space of $L^2_k$ connections
$A$ on the $\SO(3) = \PU(2)$ bundle $\su(E)$.  We shall often pass back and
forth between these viewpoints, via the fixed connection on $\det E$,
relying on the context to make the distinction clear.  Let
$D_A:L^2_k(W^+\otimes E)\to L^2_{k-1}(W^-\otimes E)$ be the corresponding
Dirac operators.  Given a connection $A$ on $E$ with curvature
$F_A\in L^2_{k-1}(\La^2\otimes\fu(E))$, then $(F_A^+)_0 \in
L^2_{k-1}(\La^+\otimes\su(E))$ denotes the traceless part of its self-dual
component. Equivalently, if $A$ is a connection on $\su(E)$ with
curvature $F_A\in L^2_{k-1}(\La^2\otimes\so(\su(E)))$, then
$\ad^{-1}(F_A^+) \in L^2_{k-1}(\La^+\otimes\su(E))$ is its self-dual
component, viewed as a section of $\La^+\otimes\su(E)$ via the isomorphism
$\ad:\su(E)\to\so(\su(E))$.

For an $L^2_k$ section $\Phi$ of $W^+\otimes E$, let $\Phi^*$ be its
pointwise Hermitian dual and let $(\Phi\otimes\Phi^*)_{00}$ be the
component of the Hermitian endomorphism $\Phi\otimes\Phi^*$ of $W^+\otimes
E$ which lies in $\su(W^+)\otimes\su(E)$. The \spinc structure $\rho$
defines an isomorphism $\rho^+:\La^+\to\su(W^+)$ and thus an isomorphism
$\rho^+=\rho^+\otimes\id_{\su(E)}$ of $\La^+\otimes\su(E)$ with
$\su(W^+)\otimes\su(E)$. Then
\begin{align}
(F_A^+)_0 - (\rho^+)^{-1}(\Phi\otimes\Phi^*)_{00} &= 0, 
\label{eq:UnpertPT}\\
D_A\Phi &= 0, \notag
\end{align}
are essentially the unperturbed equations considered in
\cite{PTCambridge,PTLocal,OTVortex,OTQuaternion} for a pair $(A,\Phi)$
consisting of a fixed-determinant connection $A$ on $E$ and a section
$\Phi$ of $W^+\otimes E$.  (The trace conditions and precise setting vary;
the equations \eqref{eq:UnpertPT} are closer to those of
\cite{TelemanNonabelian,TelemanMonopole} than \cite{PTLocal}.)
Equivalently, given a pair $(A,\Phi)$ with $A$ a connection on $\su(E)$,
the equations
\eqref{eq:UnpertPT} take the same form except that $(F_A^+)_0$ is replaced
by $\ad^{-1}(F_A^+)$ or simply by $F_A^+$, with the isomorphism
$\ad:\su(E)\to\so(\su(E))$ being implicit.

In this section we briefly describe the 
holonomy perturbations of these equations which we introduced in
\cite{FL1}: these perturbations allow us to prove transversality for the
moduli space of solutions, away from points where the connection is
reducible or the spinor vanishes identically,
and to prove the existence of an Uhlenbeck
compactification for this perturbed moduli space.

Donaldson's proof of the connected sum theorem for his polynomial
invariants \cite[Theorem B]{DonPoly} makes use of certain `extended
anti-self-dual equations' \cite[Equation (4.24)]{DonPoly}
to which the Freed-Uhlenbeck generic metrics
theorem does not apply \cite[\S 4(v)]{DonPoly}. 
To obtain transversality for the zero locus of
these extended equations, he employs holonomy perturbations which give
gauge-equivariant $C^\8$ maps $\sA_E^*\to
L^2_{k-1}(\Lambda^+\otimes\su(E))$ \cite[\S 2]{DonOrient},
\cite[pp. 282--287]{DonPoly}.  
These perturbations are continuous across the
Uhlenbeck boundary and yield transversality not only for the top stratum,
but also for all lower strata and for all intersections of the geometric
representatives defining the Donaldson invariants. 

In \cite{FL1} we describe a generalization of
Donaldson's idea which we use to prove transversality for the moduli space
of solutions to a perturbed version of the $\PU(2)$ monopole equations
\eqref{eq:UnpertPT}.
Unfortunately, in the case of the moduli space of $\PU(2)$ monopoles, the
analysis is considerably more intricate. In Donaldson's application, some
important features ensure that the requisite analysis is relatively
tractable: (i) reducible connections can be excluded from the
compactification of the extended moduli spaces \cite[p. 283]{DonPoly}, (ii)
the cohomology groups for the elliptic complex of his extended equations
have simple weak semi-continuity properties with respect to Uhlenbeck
limits
\cite[Proposition 4.33]{DonPoly}, and (iii) the perturbed zero-locus 
is cut out of a finite-dimensional manifold \cite[p. 281, Lemma 4.35, \&
Corollary 4.38]{DonPoly}. For the development of Donaldson's method for
$\PU(2)$ monopoles described here and in detail in \cite{FL1}, none of
these simplifying features hold and so the corresponding transversality
argument is rather complicated. Indeed, one can see from
Proposition 7.1.32 in
\cite{DK} that because of the Dirac operator, the behavior of the cokernels
of the linearization of the $\PU(2)$ monopole equations can be quite
involved under Uhlenbeck limits. The method we describe below uses an
infinite sequence of perturbing sections defined on the
infinite-dimensional configuration space of pairs; when restricted to small
enough open balls in the configuration space, away from reducibles, only
finitely many of these perturbing sections are non-zero and they vanish
along the reducibles.

We shall describe these perturbations and their properties only in fairly
general terms here, as the full description is lengthy; we refer the
interested reader to \cite{FL1} for a detailed account.  

Let $\sG_E$ be the Hilbert Lie group of $L^2_{k+1}$ unitary gauge
transformations of $E$ with {\em determinant one\/}. It is generally
convenient to take quotients by a slightly larger symmetry group than
$\sG_E$ when discussing pairs, so let $S_Z^1$ denote the center of $\U(2)$
and set
$$
\ssG_E := S_Z^1\times_{\{\pm\id_E\}}\sG_E,
$$ 
which we may view as the group of $L^2_{k+1}$ unitary gauge transformations
of $E$ with {\em constant determinant\/}.  The stabilizer of a unitary
connection on $E$ in $\ssG_E$ always contains the center
$S^1_Z\subset\U(2)$. We call $A$ {\em irreducible\/} if its stabilizer is
exactly $S_Z^1$ and {\em reducible\/} otherwise. Let $\sB_E(X) =
\sA_E(X)/\sG_E$ be the quotient space of $L^2_k$
connections on $E$ with fixed-determinant connection and let $\sA^*_E(X)$
and $\sB_E^*(X)$ be the subspace space of irreducible $L^2_k$ connections
and its quotient. As before, we may equivalently view $\sB_E(X)$ and
$\sB_E^*(X)$ as quotients of the spaces of $L^2_k$ connections on $\su(E)$
by the induced action of $\sG_E$ on $\su(E)$.
 
We construct gauge-equivariant $C^\8$ maps
\begin{align}
&\sA_E(X) \ni A\mapsto \vectau\cdot\vecfm(A)
\in L^2_{k+1}(X,\gl(\La^+)\otimes_\RR\so(\su(E))), 
\label{eq:GaugeEquivariantMap}\\
&\sA_E(X) \ni A\mapsto \vecvartheta\cdot\vecfm(A)
\in L^2_{k+1}(X,\Hom(W^+,W^-)\otimes_\CC\fsl(E)), \notag
\end{align}
where $\vectau = (\tau_{j,l,\alpha})$ is a sequence in
$\Om^0(X,\gl(\La^+))$ and $\vecvartheta = (\vartheta_{j,l,\alpha})$ is a
sequence in $\Om^1(X,\CC)$, while $\vecfm(A) = (\fm_{j,l,\alpha}(A))$ is
a sequence in $L^2_{k+1}(X,\su(E))$ of holonomy sections constructed by
extending the method of
\cite{DonOrient,DonPoly}, and
\begin{align*}
\vectau\cdot\vecfm(A)
&:= 
\sum_{j,l,\alpha}\tau_{j,l,\alpha}\otimes_\RR
\ad(\fm_{j,l,\alpha}(A)), \\
\vecvartheta\cdot\vecfm(A)
&:=
\sum_{j,l,\alpha}\rho(\vartheta_{j,l,\alpha})
\otimes_\CC\fm_{j,l,\alpha}(A).
\end{align*}
To construct these maps,
we fix a collection of $N_b$ small,
disjoint balls $\{4B_j\}_{j=1}^{N_b}$ in $X$, a locally finite open cover
$\{U_{j,\alpha}\}_{\alpha=1}^\8$ of each quotient space
$\sB_E^*(2B_j)$ of irreducible connections over $2B_j$, and three loops
$\{\gamma_{j,l,\alpha}\}_{l=1}^3\subset 2B_j$ such that holonomy
around these loops spans $\su(E)|_{B_j}$ for each connection in
$\{U_{j,\alpha}\}$. The sections $\fm_{j,l,\alpha}$ are supported on
$\barB_j$ in $X$ and on $L^2_k$ balls containing 
$U_{j,\alpha}$ in $\sB_E^*(2B_j)$, by a
suitable choice of cutoff functions on $X$ and $\sB_E^*(2B_j)$. The set
$\{\fm_{j,l,\alpha}(A)\}_{l=1}^3$ spans $\su(E)|_{B_j}$ for each
point $[A|_{2B_j}]\in U_{j,\alpha}$ with energy
$\|F_A\|_{L^2(4B_j)}^2<\half\eps_0^2$, where 
$\eps_0$ is a certain universal constant
\cite{FL1}. When this (regularized) energy bound is exceeded over a
ball $4B_{j'}$, the associated perturbations vanish, ensuring continuity
across the Uhlenbeck boundary.  The number $N_b$ of balls $B_j$ may be
chosen sufficiently large that for every solution $(A,\Phi)$ to the
perturbed $\PU(2)$ monopole equations
\eqref{eq:PT}, there is at least one ball $B_{j'}$
whose associated holonomy sections $\{\fm_{j',l,\alpha}(A)\}_{l=1}^3$
span $\su(E)|_{B_{j'}}$. We use the small-time heat kernel for the Neumann
Laplacians $d_A^*d_A$ on $L^2(2B_j,\su(E))$ to ensure that the sections 
$\fm_{j,l,\alpha}(A)$ are in $L^2_{k+1}$ when $A|_{2B_j}$ is in $L^2_k$.

By construction, the maps $\vectau\cdot\vecfm$ and
$\vecvartheta\cdot\vecfm$ of \eqref{eq:GaugeEquivariantMap} are uniformly
$C^s$-bounded over $\sA_E^*(X)$, when $\sA_E^*(X)$ is endowed with its
$L^2_k$ metric, {\em provided $k\ge 3$} and which we shall therefore assume
for the remainder of the article.  Moreover, they are continuous with
respect to Uhlenbeck limits, just as are those of \cite{DonPoly}. Suppose
$\{A_\beta\}$ is a sequence in $\sA_E(X)$ which converges to an Uhlenbeck
limit $(A,\bx)$ in $\sA_{E_{-\ell}}(X)\times\Sym^\ell(X)$, where
$E_{-\ell}$ is a Hermitian two-plane bundle over $X$ such that
$$
\det(E_{-\ell}) = \det E
\quad\text{and}\quad 
c_2(E_{-\ell}) = c_2(E)-\ell, 
\quad\text{with}\quad 
\ell\in\ZZ_{\ge 0}.
$$
The sections
$\vectau\cdot\vecfm(A_\beta)$ and $\vecvartheta\cdot\vecfm(A_\beta)$ then
converge in $L^2_{k+1}(X)$ to a section $\vectau\cdot\vecfm(A,\bx)$ of
$\gl(\La^+)\otimes\so(\su(E_{-\ell}))$ and a section
$\vecvartheta\cdot\vecfm(A,\bx)$ of $\Hom(W^+,W^-)\otimes\fsl(E_{-\ell})$,
respectively.  For each $\ell\ge 0$, the maps of
\eqref{eq:GaugeEquivariantMap} extend continuously to gauge-equivariant
maps
\begin{align}
&\sA_{E_{-\ell}}(X)\times\Sym^\ell(X) 
\to L^2_{k+1}(X,\gl(\La^+)\otimes_\RR\so(\su(E_{-\ell}))), 
\label{eq:GaugeEquivariantExtendedMap}\\
&\sA_{E_{-\ell}}(X)\times\Sym^\ell(X) 
\to L^2_{k+1}(X,\Hom(W^+,W^-)\otimes_\CC\fsl(E_{-\ell})), \notag
\end{align}
given by $(A,\bx)\mapsto \vectau\cdot\vecfm(A,\bx)$ and 
$(A,\bx)\mapsto \vecvartheta\cdot\vecfm(A,\bx)$, respectively,
which are $C^\8$ on each $C^\8$ stratum determined by $\Sym^\ell(X)$.

The parameters $\vectau$ and $\vecvartheta$ vary in the Banach spaces of
$\ell^1_\delta(\AAA)$ sequences in $C^r(X,\gl(\La^+))$ and
$C^r(\La^1\otimes\CC)$, 
respectively, where $\AAA = \{(j,l,\alpha)\}$ and $r\ge k+1$,
$$
\|\vecvartheta\|_{\ell^1_\delta(C^r(X))}
:= 
\sum_{j,l,\alpha}\delta_\alpha^{-1}
\|\vartheta_{j,l,\alpha}\|_{C^r(X)},
$$
and similarly for $\|\vectau\|_{\ell^1_\delta(C^r(X))}$.
The sequence of weights
$\delta=(\delta_\alpha)_{\alpha=1}^\8\in\ell^\8((0,1])$ may be chosen so
that the gauge-equivariant maps of \eqref{eq:GaugeEquivariantMap} are
smooth even at reducible connections, where the maps vanish
\cite{FL1}. 

We call an $L^2_k$ pair $(A,\Phi)$ in the {\em pre-configuration space\/},
$$
\tsC_{W,E} := \sA_E\times L^2_k(X,W^+\otimes E),
$$
a {\em $\PU(2)$ monopole\/} if it solves
\begin{align}
(F_A^+)_0 - (\id+\tau_0\otimes\id_{\su(E)} + \vectau\cdot\vecfm(A))
(\rho^+)^{-1}(\Phi\otimes\Phi^*)_{00} &= 0, 
\label{eq:PT} \\
D_A\Phi + \rho(\vartheta_0)\Phi
+\vecvartheta\cdot\vecfm(A)\Phi &= 0. \notag
\end{align}
For convenience, we often denote the perturbed
Dirac operator $D_A + \rho(\vartheta_0)+\vecvartheta\cdot\vecfm(A)$ simply by
$D_{A,\vecvartheta}$.  We let
$M_{W,E}$ be the moduli space of solutions cut out of
the {\em configuration space\/},
$$
\sC_{W,E} := \tsC_{W,E}/\ssG_E,
$$ 
by the equations \eqref{eq:PT}, where $u\in\ssG_E$ acts by $u(A,\Phi) :=
(u_*A,u\Phi)$.

We let $\sC_{W,E}^{*,0}\subset \sC_{W,E}$ be the subspace of pairs
$[A,\Phi]$ such that $A$ is irreducible and the section $\Phi$ is not
identically zero and set $M_{W,E}^{*,0} = M_{W,E}\cap
\sC_{W,E}^{*,0}$. Note that we have a canonical inclusion
$\sB_E\subset \sC_{W,E}$ given by $[A]\mapsto [A,0]$ and
similarly for the pre-configuration spaces.

The sections $\vectau\cdot\vecfm(A)$ and $\vecvartheta\cdot\vecfm(A)$
vanish at reducible connections $A$ by construction; plainly, the terms in
\eqref{eq:PT} involving the perturbations $\vectau\cdot\vecfm(A)$ and
$\vecvartheta\cdot\vecfm(A)$ are zero when $\Phi$ is zero. The holonomy
perturbations considered by Donaldson in \cite{DonPoly} are inhomogeneous,
as he uses the perturbations to kill the cokernels of $d_A^+$ directly.
In contrast, the perturbations we consider in \eqref{eq:PT} are homogeneous
and we argue indirectly that the cokernels of the linearization vanish away
from the reducibles and zero-section solutions. 

A careful application of the Agmon-Nirenberg unique continuation theorem
\cite{AN} to \eqref{eq:PT} ensures that a $\PU(2)$ monopole $(A,\Phi)$
which is irreducible on $X$ gives at least one restriction
$A|_{2B_{j'}}$ which is irreducible and whose associated holonomy
sections span $\su(E)|_{B_{j'}}$. The corresponding property for
anti-self-dual connections is proved as Lemma 4.3.21 in \cite{DK}.
The proof of Donaldson and Kronheimer relies on the Agmon-Nirenberg
unique continuation theorem for an ordinary differential inequality on
a Hilbert space \cite[Theorem 2]{AN}. We show in
\cite{FL1} that Donaldson and Kronheimer's argument adapts to the case
of the $\PU(2)$ monopole equations
\eqref{eq:UnpertPT} or \eqref{eq:PT}, when the initial open set 
where $(A,\Phi)$ is reducible contains the closed balls $\barB(x_j,R_0)$
supporting holonomy perturbations.

The perturbations $(\tau_0,\vartheta_0,\vectau,\vecvartheta)$ then ensure
that an element in the cokernel of the linearization of the parametrized
version of \eqref{eq:PT}, at a point
$(A,\Phi,\tau_0,\vartheta_0,\vectau,\vecvartheta)$ where $A$ is irreducible
and $\Phi\not\equiv 0$, must vanish identically over at least one ball
$B_{j'}$ and so must vanish identically over $X$ by the Aronszajn-Cordes
unique continuation theorem \cite{Aron}. Hence, the Sard-Smale theorem
yields:

\begin{thm}
\label{thm:Transversality}
\cite{FL1}
Let $X$ be a closed, oriented, smooth four-manifold with $C^\8$
Riemannian metric, \spinc structure $(\rho,W)$ with \spinc connection,
and a Hermitian line bundle $\det E$ with unitary connection.  Then
there exists a first-category subset of the space of $C^\8$
perturbation parameters such that the following holds: For each
$4$-tuple $(\tau_0,\vartheta_0,\vectau,\vecvartheta)$ in the
complement of this first-category subset, the moduli space
$M^{*,0}_{W,E}(\tau_0,\vartheta_0,\vectau,\vecvartheta)$ is a smooth
manifold of the expected dimension
\begin{align*} 
\dim M^{*,0}_{W,E} 
&= -2p_1(\su(E))-\threehalf(e(X)+\sigma(X)) \\
&\quad + \half p_1(\su(E))+\half(F^2-\sigma(X))-1,
\end{align*}
where $p_1(\su(E)) = c_1(E)^2-4c_2(E)$ and $F := c_1(W^+)+c_1(E)$.
\end{thm}

\begin{rmk}
Different approaches to the question of transversality for the equations
\eqref{eq:UnpertPT} with generic perturbation parameters have also
been considered by the authors, by Pidstrigach and Tyurin in \cite{PTLocal}
and by Teleman in \cite{TelemanMonopole}: see \cite{FL1} for further
details.
\end{rmk}

We now turn to the question of compactness of $M_{W,E}$, for the given
generic parameters $(\tau_0,\vartheta_0,\vectau,\vecvartheta)$. We say that
a sequence of points $[A_\beta,\Phi_\beta]$ in $\sC_{W,E}$ {\em
converges\/} to a point $[A,\Phi,\bx]$ in
$\sC_{W,E_{-\ell}}\times\Sym^\ell(X)$ if the following hold:
\begin{itemize}
\item There is a sequence of determinant-one,
$L^2_{k+1,\loc}$ bundle maps $u_\beta:E|_{X\less\{\bx\}}\to
E_{-\ell}|_{X\less\{\bx\}}$ such that the sequence of monopoles
$u_\beta(A_\beta,\Phi_\beta)$ converges to
$(A,\Phi)$ in $L^2_{k,\loc}$ over $X\less\{\bx\}$, and 
\item The sequence of measures  
$|F_{A_\beta}|^2$ converges
in the weak-* topology on measures to $|F_A|^2 +
8\pi^2\sum_{x\in\bx}\delta(x)$.
\end{itemize}
We let $M_{W,E_{-\ell}}(\bx)$ denote the moduli space of pairs
$(A,\Phi)$ solving
\eqref{eq:PT} with perturbing sections $\vectau\cdot\vecfm(\cdot,\bx)$
and $\vecvartheta\cdot\vecfm(\cdot,\bx)$, let $\bM_{W,E_{-\ell}}$
denote the moduli space of triples $(A,\Phi,\bx)$ solving \eqref{eq:PT}
for $\ell\ge 0$, and let $\bM_{W,E_{-0}} = M_{W,E}$. We define
$\barM_{W,E}$ to be the Uhlenbeck closure of $M_{W,E}$ in the space of
ideal $\PU(2)$ monopoles, 
$$
IM_{W,E} := \bigcup_{\ell=0}^N \bM_{W,E_{-\ell}} 
\subset 
\bigcup_{\ell=0}^N \left(\sC_{W,E_{-\ell}}\times\Sym^\ell(X)\right)
$$
for any integer $N\ge N_p$ where $N_p$ is a sufficiently large
constant. Analogues of Bochner formulas used in the proof of compactness
for the Seiberg-Witten equations \cite{KMThom,Witten} provide a universal
energy bound for solutions to \eqref{eq:PT}, guaranteeing that the
constants $N_b$ and $N_p$ exist. By combining the methods used in the proof
of compactness for the Seiberg-Witten moduli space \cite{KMThom} and
Uhlenbeck compactness for the moduli space of anti-self-dual equations
\cite{DK} we obtain: 

\begin{thm}\label{thm:Compactness}
\cite{FL1}
Let $X$ be a closed, oriented, smooth four-manifold with $C^\8$ Riemannian
metric, \spinc structure $(\rho,W)$ with 
\spinc connection, and a Hermitian two-plane
bundle $E$ with unitary connection on $\det E$.  Then there is a positive
integer $N_p$, depending at most on the curvatures of the fixed connections
on $W$ and $\det E$ together with $c_2(E)$, such that for all $N\ge N_p$
the topological space $\barM_{W,E}$ is compact, second-countable,
Hausdorff, and is given by the closure of $M_{W,E}$ in
$\cup_{\ell=0}^{N}\bM_{W,E_{-\ell}}$.
\end{thm}

\begin{rmk}
The existence of an Uhlenbeck compactification for the moduli space of
solutions to the unperturbed $\PU(2)$ monopole equations
\eqref{eq:UnpertPT} was announced by Pidstrigach
\cite{PTCambridge} and an argument was outlined in \cite{PTLocal}.
A similar argument for the equations 
\eqref{eq:UnpertPT} was outlined by Okonek and Teleman in
\cite{OTQuaternion}. Theorem \ref{thm:Compactness} yields the standard
Uhlenbeck compactification for the system \eqref{eq:UnpertPT} and for the
perturbations of \eqref{eq:UnpertPT} described in \cite{PTLocal}.
A proof of Uhlenbeck compactness for 
\eqref{eq:UnpertPT} (and for certain perturbations of these equations) is
also given in \cite{TelemanMonopole}.
\end{rmk}

We use the term {\em (Uhlenbeck) level\/} to describe the spaces
$\bM_{W,E_{-\ell}}$ for different values of $\ell\ge 0$, with $M_{W,E}$
comprising the {\em top (Uhlenbeck) level\/}. The space $\Sym^\ell(X)$ is
smoothly stratified, the strata being enumerated by partitions of
$\ell$. If $\Sigma\subset\Sym^\ell(X)$ is a smooth stratum, we define
$$
\bM_{W,E_{-\ell}}|_\Sigma
:=\{[A,\Phi,\bx]\in\bM_{W,E_{-\ell}}: \bx\in\Sigma \}.
$$
The proof of Theorem \ref{thm:Transversality} shows, more generally, that 
for each $\ell\ge 0$ the moduli spaces
$$
\bM^{*,0}_{W,E_{-\ell}}|_\Sigma
:= \bM_{W,E_{-\ell}}|_\Sigma\cap \bM^{*,0}_{W,E_{-\ell}}
$$ 
are smooth and of the expected dimension, and over the complement in
$\Sigma$ of a first-category subset, the projection
$\bM_{W,E_{-\ell}}^{*,0}|_\Sigma\to\Sigma$ is a fiber bundle.  See
\cite{FL1} for the general statement. In the more familiar case of the
Uhlenbeck closure of the moduli space of solutions to the unperturbed
$\PU(2)$ monopole equations
\eqref{eq:UnpertPT}, the spaces $\bM_{W,E_{-\ell}}$
would be replaced by the products $M_{W,E_{-\ell}}\times \Sym^\ell(X)$.  In
general, though, the spaces $\bM_{W,E_{-\ell}}$ are not
products due to the slight dependence of the lower-level analogues of
the equations \eqref{eq:PT} on the points $\bx\in\Sym^\ell(X)$. A similar
phenomenon is encountered in \cite[\S 4(iv)--(vi)]{DonPoly} for the
case of the extended anti-self-dual equations.

While the description of the holonomy perturbations outlined above may
appear fairly complicated at first glance in practice, they do not present
any major difficulties beyond those that would be encountered if simpler
perturbations not involving the bundle $\su(E)$ (such as the Riemannian
metric on $X$ or the connection on $\det W^+$) were sufficient to achieve
transversality \cite{FL2, FL3, FL4}. We note that related transversality
and compactness issues have been recently considered in approaches to
defining Gromov-Witten invariants for general symplectic manifolds
\cite{LiTian,RuanGW,Siebert}.


\section{Cohomology and cobordisms}
\label{sec:Cohomology}
The moduli space $M_{W,E}$ contains singularities: it is a smoothly
stratified space, with strata diffeomorphic to the moduli space of
anti-self-dual connections on $\su(E)$
and to moduli spaces of $\U(1)$ monopoles (which
are in turn cobordant to moduli spaces of Seiberg-Witten monopoles). The
space $M^{*,0}_{W,E}$ therefore gives a cobordism between the links of these
two types of singularities. In this section, we introduce cohomology
classes on $M^{*,0}_{W,E}$ and define the links of these singularities.

\subsection{Singularities}
\label{sec:Singularities}
We see from Theorem \ref{thm:Transversality}
that the moduli space $M^{*,0}_{W,E}$ of $\PU(2)$ monopoles
$[A,\Phi]$, where $A$ is not reducible and $\Phi\not\equiv 0$, forms a
smooth manifold.  We now describe the subspaces where $A$ is reducible or
$\Phi\equiv 0$.  

Let $M^{\asd}_E\subset M_{W,E}$ denote the subspace of points $[A,\Phi]$
where $\Phi\equiv 0$; we refer to pairs representing points in $M^{\asd}_E$
as {\em zero-section pairs}. Equivalently, we may view $M^{\asd}_E\subset
\sB_E$ as the moduli space of fixed-determinant connections $A$ on $E$
solving the {\em anti-self-dual equation\/},
\begin{equation}
(F_A^+)_0 = 0,
\label{eq:ASD}
\end{equation}
or simply $F_A^+ = 0$, if $\sB_E$ is viewed as the quotient space of 
connections $A$ on $\su(E)$.

Suppose we have a reduction of the $\U(2)$ bundle $E$ given as an (ordered)
direct sum of line bundles,
$$
E=L_1\oplus L_2.
$$
Note that gauge transformations of $E$ (in $\ssG_E =
S_Z^1\times_{\{\pm\id_E\}}\sG_E$) which interchange the line bundles $L_1$ and
$L_2$ only exist if $L_1 = L_2$.
We let $M^{\red}_{W,E,L_1}\subset M_{W,E}$ denote the subspace of points
$[A,\Phi]$ with $\Stab_{A,\Phi}=S^1_{L_2}$, where $S^1_{L_2} = S^1$ acts by
constant multiplication on the line bundle $L_2$. We refer to pairs
representing points in $M^{\red}_{W,E,L_1}$ as {\em reducible pairs}: they
have the form $(A_1\oplus A_2,\Phi_1)$, where $A_1$ is a unitary connection
on $L_1$ and $A_2 = A_e\otimes A_1^*$ is the corresponding connection on
$L_2 = (\det E)\otimes L_1^*$, where $A_e$ is the fixed connection on $\det E$,
while $\Phi_1$ is a section of $W^+\otimes
L_1$. The pair $(A_1,\Phi_1)$ is a solution to the {\em $\U(1)$ monopole
equations\/},
\begin{align}
F_{A_1}^+ - \half(\id+\tau_0)(\Phi_1\otimes\Phi_1^*)_0 - \half F_{A_e}^+
&= 0, \label{eq:U1Monopole} 
\\
D_{A_1}\Phi_1 &=0. \notag
\end{align}
The moduli space of solutions to \eqref{eq:U1Monopole}, which parametrizes
$M^{\red}_{W,E,L_1}$, is smooth and of the expected dimension for generic
$\tau_0$ away from the zero-section solutions (see \cite{FL2}) and 
is cobordant to the standard Seiberg-Witten moduli
space $M_{W\otimes L_1}^{\sw}$ associated to
the \spinc structure $(\rho,W\otimes L_1)$ (as defined, for example, in
\cite{MorganSWNotes}). 

\begin{prop} 
\label{prop:ClassificationOfStabilizers}
\cite{FL2}
Let $X$ be a closed, oriented, smooth four-manifold
with $b^+(X)\geq 1$ and generic Riemannian metric. Suppose the pair $(A,\Phi)$
on $(E,W^+\otimes E)$ represents a point $[A,\Phi]\in M_{W,E}$ 
with non-trivial stabilizer $\Stab_{A,\Phi}$. Then one of the following,
mutually exclusive situations holds:
\begin{enumerate}
\item 
The pair $(A,\Phi)$ is a zero-section pair ($\Phi\equiv 0$)
and the connection $A$ is
irreducible. The pair $(A,0)$ has stabilizer $\Stab_{A,0}=S^1_Z$,
the connection $A$ has stabilizer $\Stab_A=S^1_Z$, and $A$ is projectively
anti-self-dual (so $(F^+_A)_0=0$).  The quotient space of zero-section pairs is
identified with the moduli space $M_E^{\asd}$ of anti-self-dual
connections on $\su(E)$.
\item 
The pair $(A,\Phi)$ is reducible and $\Phi\not\equiv 0$.
The bundle $E$ splits as $E=L_1\oplus L_2$, the pair $(A,\Phi)$ has
stabilizer $\Stab_{A,\Phi}=S^1_{L_2}$, and $A$ has stabilizer 
$\Stab_A=S^1_{L_1}\times S^1_{L_2}$.  
If $M^{\asd}_E\cap M^{\red}_{W,E,L_1}=\emptyset$, then
$M^{\red}_{W,E,L_1}$ is smoothly cobordant to
the Seiberg-Witten moduli space 
$M^{\sw}_{W\otimes L_1}$.
\item 
The pair $(A,\Phi)$ is a reducible, zero-section pair.
The connection $A$ is projectively flat (so $(F_A)_0=0$)
and $\Phi\equiv 0$.  The bundle $E$ splits as 
$E=L_1\oplus (L_1\otimes N)$, where $N$ is a torsion line bundle, so
$c_1(N)\in\Tor H^2(X;\ZZ)$.
The stabilizer of the pair is $\Stab_{A,0}=\Stab_A$.
\end{enumerate}
If $b^+(X)=0$ or the Riemannian metric
metric on $X$ is non-generic, the pair $(A,\Phi)$
can have stabilizer $\Stab_{A,\Phi}=S^1_{L_1}\times S^1_{L_2}$, where
$\Phi\equiv 0$ and $A$ is a reducible projectively anti-self-dual, but not
projectively flat connection.
\end{prop}

\begin{rmk}
If $X$ is simply-connected, then the third case only occurs when the
connection on $\su(E)$ induced by $A$ is trivial.  The stabilizer of the
pair is then $\U(2)$.
\end{rmk}

The undesirable third case in Proposition
\ref{prop:ClassificationOfStabilizers} (see \cite{FL2}) can be excluded
with the aid of a criterion due to Fintushel and Stern \cite{FSFlat}:

\begin{prop}
\label{prop:ProjFlatCriterion}
\cite{FSFlat}
If $c\in H^2(X;\ZZ)$ and $c\pmod{2}\in H^2(X;\ZZ_2)$ is not a pullback
from $H^2(K(\pi_1(X),1);\ZZ_2)$, then there are no 
$\SO(3)$ bundles $V\to X$ with $w_2(V)=c\pmod{2}$
which admit a flat connection.
\end{prop}

We can choose the class $w_2(\su(E)) = c_1(E)\pmod{2}$ so that $\su(E)$ does
not admit a flat connection using the blow-up trick of \cite{MorganMrowka}:
If $c\in H^2(X;\ZZ)$ and $e^*$ is the Poincare dual of the exceptional
class of the blow-up $\hat X := X\#\overline{\CC\PP}^2$, then $c+e^*$ does
not admit a flat $\SO(3)$ connection. As the Donaldson polynomials and
Seiberg-Witten invariants of $X$ and its blowup $\hat X$ determine each
other, no information is lost in this process \cite{FSBlowupSW,FSBlowupD}.
Therefore, assuming this third possibility does not occur, the moduli space
$M_{W,E}$ has a smooth stratification
\begin{equation}
M_{W,E}
=
M^{*,0}_{W,E} \cup M^{\asd}_E \cup M^{\red}_{W,E},
\quad\text{with}\quad
M^{\red}_{W,E} 
:= 
\bigcup_{L_1}M^{\red}_{W,E,L_1},
\label{eq:TopLevelStratification}
\end{equation}
where the union is over the finitely many line bundles $L_1\in H^2(X;\ZZ)$
for which (i) there is a topological splitting $E=L_1\oplus L_2$, where $L_2
= (\det E)\otimes L_1^*$ and recalling that $\det E$ is fixed, and (ii)
the moduli space $M^{\red}_{W,E,L_1}$ is non-empty.  One can show
directly that there are only a finite number of line bundles $L_1$ with
$M^{\red}_{W,E,L_1}$ non-empty by repeating the usual argument for the
standard Seiberg-Witten moduli spaces \cite[Theorem 5.2.4]{MorganSWNotes}.

{\em For the remainder of this article\/} we shall assume that $X$ is
equipped with an orientation, a homology orientation, has $b^+(X)> 0$, and
is equipped with a generic Riemannian metric. In the case $b^+(X) = 1$, the
Donaldson invariants refer to the specific chamber in $H^2(X;\RR)/\RR^*$
defined by the choice of metric. The dimensions of our moduli spaces are
then given by
\begin{align*}
2d(\su(E),F) 
:= \dim M^{*,0}_{W,E} 
&= 2d_a(\su(E)) + 2n_a(\su(E),F) - 1,
\\
2d_a(\su(E))
:= \dim M^{\asd}_E
&= -2p_1(\su(E))-\threehalf(e(X)+\sigma(X))
\\
&= -2p_1(\su(E))-3(1-b^1(X)+b^+(X)),
\\
2n_a(\su(E),F) 
:= 2\Ind_\CC D_A
&= \half p_1(\su(E))+\half(F^2-\sigma(X)),
\end{align*}
where $p_1(\su(E)) = c_1(E)^2-4c_2(E)$ and $F = c_1(W^+)+c_1(E)$, while
\begin{align*}
2d_s(K)
:= \dim M^{\sw}_{W\otimes L_1} 
&= \quarter(K^2-(2e(X)+3\sigma(X))) 
\\
&= \quarter(K^2-\sigma(X)) - (1-b^1(X)+b^+(X)), 
\end{align*}
where $K := c_1(W^+\otimes L_1)$.  Since $b^+(X)-b^1(X)$ is odd for $X$
admissible, one sees that the moduli spaces $M^{\asd}_E$ and
$M^{\sw}_{W\otimes L_1}$ are indeed even-dimensional, as implied above,
since $\quarter(K^2-\sigma(X))$ is twice the complex index of the Dirac
operator on $W^+\otimes L_1$.

\subsection{Cobordisms of links via moduli spaces of PU(2) monopoles}
The essential idea is to use the moduli space
$M^{*,0}_{W,E}$ as a cobordism between the `links' of 
$M^{\asd}_E$ and $M^{\red}_{W,E}$.  In \S
\ref{subsec:CohomologyDef} we define cohomology classes and their dual
geometric representatives on
$M^{*,0}_{W,E}$.  The pairing of a product of these cohomology classes 
(or intersection of their dual geometric representatives) with
the link of $M^{\asd}_E$ can be expressed as a multiple of the 
Donaldson polynomial (Lemma \ref{lem:LinkOfASD}) while the pairing
of these classes with the link of $M^{\red}_{W,E}$ gives
multiples of the Seiberg-Witten invariants
(Theorem \ref{thm:DegreeZeroFormula}).  
The intersection of the geometric representatives in $M^{*,0}_{W,E}$
is a family of oriented one-manifolds, whose boundaries should
lie in the links of $M^{\asd}_E$ and $M^{\red}_{W,E}$,
yielding an equality between
these pairings and thus a relationship between the Donaldson and Seiberg-Witten
invariants.

Two technical difficulties arise in the above program.  
The first problem is that $M^{*,0}_{W,E}$ is not compact.
Thus the boundaries of the one-manifolds might not lie on
these links, but in the lower levels of $\barM_{W,E}$.
One can work instead with
$\barM^{*,0}_{W,E}$, 
the subspace of $\barM_{W,E}$
given by triples $[A,\Phi,\bx]$ where $\Phi\not\equiv 0$ and
$A$ is not reducible.  In \S \ref{subsec:ExtendGeomRepr},
we describe the intersection of the closure of the geometric
representatives in $\barM_{W,E}$
with the lower strata of $\barM^{*,0}_{W,E}$.
This description and a dimension-counting argument show 
that the one-manifolds given by the intersection of the
geometric representatives do not have boundary points
in the lower levels of $\barM^{*,0}_{W,E}$.

The second problem is to define links of the singularities $M^{\asd}_E$
and $M^{\red}_{W,E,L_1}$.  The equations \eqref{eq:PT} cutting out
$M_{W,E}\subset\sC_{W,E}$ do not vanish transversely along these
singularities and so the local topology of $M_{W,E}$ could be quite
intricate near $M^{\asd}_E$ and $M^{\red}_{W,E,L_1}$.  In \S
\ref{subsec:CohomOnASDLink} we define a smoothly-stratified,
codimension-one subspace $\bL^{\asd}_\varepsilon\subset
\barM^{*,0}_{W,E}$ and in 
Lemma \ref{lem:LinkOfASD} we compute the intersection of some geometric
representatives with this link. In \S \ref{subsec:LinkOfReducible} we
outline our definition \cite{FL2} of a link $\bL_{W,E,L_1}\subset
M^{*,0}_{W,E}$ of the stratum $M^{\red}_{W,E,L_1}$ in $M_{W,E}$ and
describe the intersection of the geometric representatives with this link
in Theorem
\ref{thm:DegreeZeroFormula}.

If all the reducibles lie only in the top level of $\barM_{W,E}$, the
cobordism $\barM^{*,0}_{W,E}$ yields an explicit formula relating the
Donaldson polynomial and the Seiberg-Witten invariant (Theorem
\ref{thm:CompactReductionFormula}).  In general, however, there will be
reducible pairs in the lower levels of $\barM_{W,E}$.  The one-manifolds
given by the intersection of the geometric representatives can then have
boundaries at reducible pairs in the lower levels of $\barM_{W,E}$.  The
space $\barM^{*,0}_{W,E}$ yields a cobordism between
$\bL^{\asd}_\varepsilon$ and the links of all the reducibles, including
these lower-level reducibles.  The definition of the links of the
lower-level reducibles is considerably more involved and is discussed in \S
\ref{sec:GluingKMConjecture}.

\subsection{The cohomology classes}
\label{subsec:CohomologyDef}
In this subsection we define the cohomology classes on $M^{*,0}_{W,E}$,
referring the reader to \cite{FL2} for detailed description of their
dual geometric representatives. Recall that
$\tsC_{W,E}=\sA_E\times\Omega^0(W^+\otimes E)$ is our pre-configuration
space of $L^2_k$ pairs,
where we have omitted Sobolev indices as these play
no role in the present discussion.  Let $\tsC^*_{W,E}$
denote the subspace of pairs which are not reducible, let $\tsC^0_{W,E}$
denote the subspace of those which are not zero-section pairs, and let
$\tsC^{*,0}_{W,E}$ denote the subspace of those
which are neither zero-section nor reducible pairs.
Let $P$ be the $\U(2)$ principal bundle underlying the vector bundle $E$
and define
$$
\PP := \tsC^{*,0}_{W,E}\times_{\ssG_E}P.
$$ 
The space $\PP$ is a principal $\U(2)$ bundle over 
$\sC^{*,0}_{W,E}\times X$.  
The associated $\SO(3)$ bundle,
$\PP^{\ad}:=\PP/S^1_Z$, extends over $\sC^*_{W,E}$.
Indeed, the space $\PP$
is isomorphic to $\PP/S^1_Z$ over the zero-section pairs.  
Over the reducible pairs, 
the space $\PP$ becomes an $\SO(3)$ fiber bundle, 
but is not principal as the 
stabilizers of these pairs are not normal subgroups of $\U(2)$.

We define maps from the homology of $X$ to the cohomology of
$\sC^{*,0}_{W,E}$ via
\begin{align*}
&\mu_{c_1}: H_\bullet(X;\QQ)\rightarrow H^{2-\bullet}(\sC^{*,0}_{W,E};\QQ),
\quad \beta \mapsto c_1(\PP)/\beta, \\
&\mu_{p_1}: H_\bullet(X;\QQ)\rightarrow H^{4-\bullet}(\sC^{*}_{W,E};\QQ), 
\quad \beta\mapsto -\quarter p_1(\PP/S^1_Z)/\beta,
\end{align*}
where 
$$
\quarter p_1(\PP/S^1_Z)/\beta
                 = (c_2(\PP)-\quarter c_1^2(\PP))/\beta.
$$
Following \cite[Definition 5.1.11]{DK} we define a universal $\SO(3)$
bundle by
$$
\PP_E^{\ad} 
:= 
\sA_E^*\times_{\sG_E}(P/S^1_Z)
\rightarrow 
\sB^*_E\times X
$$
and set
$$
\mu_E: H_\bullet(X;\QQ)\rightarrow H^{4-\bullet}(\sB^*_E;\QQ),
\quad \beta \mapsto -\quarter p_1(\PP_E^{\ad})/\beta.
$$
If $\pi:\sC^{*,0}_{W,E}\rightarrow \sB^*_E$ is the projection
$[A,\Phi] \mapsto [A]$, we see that $(\pi\times
\id_X)^*\PP_E^{\ad} = \PP^{\ad}$. This implies the following relation
between the cohomology classes on $\sC^{*,0}_{W,E}$ and $\sB^*_E$:

\begin{lem}
\label{lem:PullbackUniversalSO3Bundle}
If $\beta\in H_\bullet(X;\QQ)$, then $\pi^*\mu_E(\beta)=\mu_{p_1}(\beta)$.
\end{lem}

The class $\mu_{c_1}(x)$ is non-trivial on the link of the zero-section
pairs \cite{FL2}.  It does not pull back from the quotient
space of connections and does not even extend over the subspace
$M_E^{\asd}\subset M_{W,E}$.

By analogy with the construction of geometric representatives for
cohomology classes in Donaldson theory
\cite{DonConn,DonPoly,DK,KMStructure}, we define geometric representatives
$V(\beta)$ and $W(x)$ to represent $\mu_{p_1}(\beta)$ and $\mu_{c_1}(x)$,
respectively.  Some features of the definition of these geometric
representatives are worth mentioning.  For a smooth submanifold $Y\subset
X$ representing $\beta\in H_\bullet(X;\QQ)$, we let $U_Y$ be a `suitable'
neighborhood \cite[\S 2]{KMStructure}.  The representatives $V(\beta)$ are
the pull-backs of the usual usual geometric representatives of Donaldson
theory \cite{KMStructure} from the quotient space of connections
$\sB_E^*(U_Y\cup_j B_j)$, where $\barB_j$ are the balls supporting the
holonomy perturbations.  If the energy of a connection $A|_{4B_{j'}}$ is
greater than a certain universal bound, the representative $V(\beta)$ is
independent of its restriction to $B_{j'}$.

As in \cite{KMStructure}, we let
$\AAA(X) := \Sym\left( H_{\even}(X;\QQ)\right)\otimes
\La\left( H_{\odd}(X;\QQ)\right)$
be the graded algebra, with $z=\beta_1\beta_2\dots\beta_r$ having total
degree $\deg(z) = \sum_i(4-i_p)$,
when $\beta_p\in H_{i_p}(X;\QQ)$.  We write
\begin{align*}
\mu_{p_1}(z) &:= \mu_{p_1}(\beta_1)\smile\dots\smile\mu_{p_1}(\beta_r), \\
\\
V(z) &:= V(\beta_1)\cap\dots\cap V(\beta_r),
\end{align*}
for $z=\beta_1\beta_2\dots\beta_r$, and similarly for $\mu_E(z)$. We write
$$
\mu_{c_1}(x^m) 
:= 
\underbrace{\mu_{c_1}(x)\smile\dots\smile\mu_{c_1}(x)}_{\text{$m$ times}} 
\quad\text{and}\quad
W(x^m) := \underbrace{W(x)\cap\cdots\cap W(x)}_{\text{$m$ times}},
$$
for products of the class $\mu_{c_1}(x)$ and its dual $W(x)$. 

\subsection{The closure of the geometric representatives}
\label{subsec:ExtendGeomRepr}
We now describe the intersection of the geometric representatives
with the lower strata of $\barM_{W,E}$.
Let $\Sigma\subset\Sym^\ell(X)$ be a smooth stratum.
Counting dimensions, one sees that 
\begin{align*}
\dim \bM^{*,0}_{W,E_{-\ell}}(\Sigma)
&= \dim M^{*,0}_{W,E} - 6\ell+\dim\Sigma \\
&\leq 
\dim M^{*,0}_{W,E}-2\ell, \qquad 0 \le \ell \le N_p,
\end{align*}
so the strata $\bM^{*,0}_{W,E_{-\ell}}(\Sigma)$ (with $\ell\ge 1$) of the
compactification $\barM_{W,E}$ have codimension at least two less
than the top stratum $M^{*,0}_{W,E}$.  This would allow the
definition of a relative fundamental class (with boundaries given by the
links of the zero-section and reducible pairs) if we knew
$\barM_{W,E}$ had locally finite topology.  We consider
intersections of geometric representatives whose total codimension is one
less than the dimension of $M^{*,0}_{W,E}$.  Thus, if these
geometric representatives intersect the lower strata of
$\barM_{W,E}$ in sets of the same codimension as their intersection
with the top stratum $M^{*,0}_{W,E}$, dimension counting shows that
the intersection of these geometric representatives, away from the
zero-section and reducible pairs, occurs only in the top stratum.

\begin{defn}
The closures of the geometric representatives, $V(\beta)$, $W(x)$, in
$\barM_{W,E}$ are denoted by $\barV(\beta)$, $\barW(x)$, respectively.  For
$z=\beta_1\dots\beta_r\in\AAA(X)$, a generator $x\in H_0(X)$, and an
integer $m\ge 0$, we denote
$$
\barV(z) := \barV(\beta_1)\cap\cdots\cap\barV(\beta_r)
\quad\text{and}\quad
\barW(x^m) := \underbrace{\barW(x)\cap\cdots\cap\barW(x)}_{\text{$m$ times}}.
$$
\end{defn}
  
The description of the intersection of $\barV(\beta)$, $\barW(x)$ with the
lower strata given below in Lemma \ref{lem:CyclesExtension} is incomplete,
as it (i) gives only an inclusion and not an equality and (ii) does not
give the multiplicities of components of these intersections occuring in
lower levels.  A more complete description is given in \cite{FL4}, using
`tubular neighborhood' descriptions of the lower strata in $\barM_{W,E}$
obtained from gluing maps.  

For $i=1,\dots,\ell$, let $\pi_i:X\times\dots\times X\rightarrow X$ be
projection onto the $i$th factor.  Let $S^\ell(Y)$ be the projection of
$\cup_i\pi_i^{-1}(Y)$ to $\Sym^\ell(X)$ under the map $X^\ell\rightarrow
\Sym^\ell(Y)$ and denote $S_\Sigma(Y)=\Sym^\ell(Y)\cap \Sigma$.

On each space $\bM_{W,E_{-\ell}}^{*,0}$, there are geometric
representatives $V_\ell(\beta)$ and $W_\ell(x)$ defined in exactly the same
way as the geometric representatives $V(\beta)$, $W(x)$ on $M^{*,0}_{W,E}$,
except that we use bundles $P_{-\ell}$ and $P_{-\ell}^{\ad} :=
(P_{-\ell})/S^1_Z$ with $c_1(P_{-\ell})=c_1(P)$ and
$c_2(P_{-\ell})=c_2(P)-\ell$.  We then have the following description of
the intersection of the extended geometric representatives $\barV(\beta)$,
$\barW(x)$ with $\bM^{*,0}_{W,E_{-\ell}}(\Sigma)$:

\begin{lem}
\label{lem:CyclesExtension}
For a smooth stratum $\Sigma\subset \Sym^\ell(X)$, let $\pi:
\bM_{W,E_{-\ell}}^{*,0} (\Sigma)\rightarrow \Sigma$ be the projection
map.  Let $x\in H_0(X)$ be a generator, let
$\beta\in H_\bullet(X;\QQ)$ have a smooth representative
$Y\subset X$, and let $U_Y$ be a suitable neighborhood of $Y$.  Then the
following hold:
\begin{enumerate}
\item 
   $\barV(\beta)\cap \bM_{W,E_{-\ell}}^{*,0} (\Sigma)
        \subseteq V_\ell(\beta)\cup \pi^{-1}(S_\Sigma(U_Y))$,
\item 
   $\barW(x)\cap \bM_{W,E_{-\ell}}^{*,0} (\Sigma)
        \subseteq W_\ell(x)\cup \pi^{-1}(S_\Sigma(U_x))$.
\end{enumerate}
Furthermore, if $\ell=0$ and $\beta\in H_2(X;\QQ)$ 
is a two-dimensional class with $\langle  2L_1-c_1(E),\beta\rangle\neq 0$,
then we have the following reverse inclusions:
\begin{enumerate}
\item 
   $M^{\red}_{W,E,L_1}\subset \barV(\beta)$,
\item 
   $M_{W,E}^{\red}\subset\barV(x)$,
\item 
   $\barM_E^{\asd}\cup M_{W,E}^{\red}\subset\barW(x)$.
\end{enumerate}
\end{lem}

\begin{rmk}
\begin{enumerate}
\item The intersections of the geometric representatives with 
the strata of reducible pairs and of zero-section pairs in $\barM_{W,E}$
generally do not have the expected codimensions.  Indeed, Lemma
\ref{lem:CyclesExtension} shows that almost all geometric
representatives will contain reducible pairs in the top level. 
\item To get equality in the first assertions
(replacing $S_\Sigma(U_Y)$ with $S_\Sigma(Y)$),
we use gluing to describe the geometric representatives in
an Uhlenbeck neighborhood of the lower level.
\end{enumerate}
\end{rmk} 

One cannot use dimension counting directly at this point as
the open subsets $\pi^{-1}(S_\Sigma(U_Y))$ 
in $\bM_{W,E_{-\ell}}^{*,0} (\Sigma)$
do not have positive codimension.
However, it can be shown that the restrictions of the
geometric representatives $V_\ell(\beta)$, $W_\ell(x)$ to
$\pi^{-1}(S_\Sigma(U_Y))$ are given by a pullback from 
$\pi^{-1}(S_\Sigma(Y))$.
The intersection of the geometric representatives
with $\bM_{W,E_{-\ell}}^{*,0}(\Sigma)$ may thus be computed by
replacing
$\pi^{-1}(S_\Sigma(U_Y))$ with $\pi^{-1}(S_\Sigma(Y))$.

We then see from Lemma \ref{lem:CyclesExtension} and the transversality
results of \S \ref{sec:TransvCompact} that although the closures
$\barV(\beta)$ and $\barW(x)$ do not intersect every stratum of
$\barM_{W,E}$ in a set of the expected codimension, they do intersect the
strata of $\barM^{*,0}_{W,E}$ in sets of the expected codimension.  A
dimension-counting argument then yields:

\begin{cor}
\label{cor:IntersectionOfSmoothLowerStrata}
\cite{FL2}
Let $n_{p_1}$ and $n_{c_1}$ be non-negative
integers such that $n_{p_1}+n_{c_1}=d_a+n_a-1$.
Let $\beta_1,\dots,\beta_r\in H_\bullet(X;\QQ)$
be homology classes such that $\sum_i (4-\dim \beta_i)=n_{p_1}$
and let $z=\beta_1\beta_2\dots\beta_r\in\AAA(X)$.
If the collection $\beta_1,\dots,\beta_r$ does not contain
both a zero-dimensional class and a three-dimensional class,
then for generic choices of geometric representatives, and
appropriate choices of suitable neighborhoods, the intersection
$$
\barV(z) \cap \barW(x^{n_{c_1}}) \cap \barM^{*,0}_{W,E}
$$
is a collection of one-dimensional manifolds, 
disjoint from the lower strata of $\barM^{*,0}_{W,E}$.
\end{cor}

\begin{rmk}
The condition in Corollary \ref{cor:IntersectionOfSmoothLowerStrata}
about the absence of either three- or zero-dimensional homology
classes is necessary because the definition of a suitable neighborhood
includes loops which weaken the conclusions one can reach by 
dimension counting (see \cite[p. 593]{KMStructure} or \cite{FL2}).
\end{rmk} 

\subsection{Orientations and the deformation complex}
\label{subsec:Orient}
The deformation complex for the $\PU(2)$ monopole equations \eqref{eq:PT}
is given by
\begin{equation}
\label{eq:DeformationComplex}
\begin{CD}
\Om^0(\su(E))\oplus i\RR_Z
@>{d_{A,\Phi}^0}>> 
\begin{matrix}
\Om^1(\su(E)) \\
\oplus \\
\Om^0(W^+\otimes E)
\end{matrix}
@>{d_{A,\Phi}^1}>> 
\begin{matrix}\Om^+(\su(E)) \\
\oplus \\
\Om^0(W^-\otimes E)
\end{matrix}
\end{CD} 
\end{equation}
where $i\RR_Z$ is the Lie algebra of $S^1_Z$.  Here, $d^0_{A,\Phi}$ is the
differential of the action of the gauge group $\ssG_E$
at $(A,\Phi)$, while $d^1_{A,\Phi}$
is the linearization of the $\PU(2)$ monopole equations \eqref{eq:PT}.  Let
$$
\sD_{A,\Phi} := d^{0,*}_{A,\Phi}+d^1_{A,\Phi}
$$ 
be the `rolled-up' deformation operator.  For any point $[A,\Phi]\in
M^{*,0}_{W,E}$, there is an isomorphism, $T_{A,\Phi}M^{*,0}_{W,E}\simeq
\Ker\sD_{A,\Phi}$. In \cite{FL2} we prove that $M^{*,0}_{W,E}$ is
orientable by showing that the real line bundle $\det\sD$ is trivial.

An orientation for $M^{*,0}_{W,E}$ can be specified by choosing a
value for a section of $\det\sD$ at any point $[A,\Phi]\in\sC_{W,E}$.
At a zero-section $\PU(2)$ monopole $(A,0)$, the deformation complex
\eqref{eq:DeformationComplex} splits into the direct sum of complexes:
\begin{gather*}
\begin{CD}
\Omega^0(\su(E))
@> d_A>>
\Omega^1(\su(E))
@> d^+_A >>
\Omega^+(\su(E)) 
\end{CD} 
\\
\begin{CD}
\Omega^0(W^+\otimes E)
@> D_A >>
\Omega^0(W^-\otimes E)
\end{CD}
\end{gather*} 
The first complex is the elliptic deformation complex for the moduli space
$M_E^{\asd}$ of anti-self-dual connections
and $i\RR_Z$ is in the cokernel of $\sD_{A,0}$. Because
\begin{equation}
\label{eq:ASDOrientIsom}
\det\sD_{A,0}\simeq 
\det\left( d^*_A+d^+_A\right)\otimes\det D_A\otimes (i\RR_Z)^*,
\end{equation}
we can specify an orientation for $\det\sD$ by specifying one
for the anti-self-dual moduli space, using the complex orientation
on $\det D_A$, and fixing an orientation for $i\RR_Z$.

\begin{defn}
If $w\in H^2(X;\ZZ)$ is an integral lift of $w_2(\su(E))$ and
$p_1(\su(E))=-4k$, and $\Omega$ is a homology orientation for $X$, let
$o_k(\Omega,w)$ be the corresponding orientation defined in \cite[\S
7.1.6]{DK} for the moduli space $M_E^{\asd}$ of anti-self-dual connections
on $\su(E)$.  Let $O^{\asd}_{k}(\Omega,w)$ be the orientation for
$\det\sD$, and so $M^{*,0}_{W,E}$, defined through the isomorphism
\eqref{eq:ASDOrientIsom}, the orientation $o_{k}(\Omega,w)$
for the moduli space $M_E^{\asd}$, the complex orientation for $\det D$ and
the fixed orientation for $i\RR_Z$. The moduli space $M^{\asd}_E$ is
equipped with the {\em standard orientation\/} $o_{k}(\Omega,c_1(E))$,
where $k=-\quarter p_1(\su(E))$, if no other orientation is specified.
\end{defn}

\begin{rmk}
Since $p_1(\su(E))=c_1(E)^2-4c_2(E)$ and $w_2(\su(E)) = c_1(E)\pmod{2}$,
then $p_1(\su(E))=w^2\pmod{4}$ if $w$ is an integral lift of $w_2(\su(E))$.
The orientation for $M_E^{\asd}$ is then determined by the addition of
$-\quarter(p_1(\su(E))-w^2)$ instantons to the $\U(2)$ bundle 
$\underline\CC\oplus w$, with corresponding $\SO(3)$ bundle
$\underline\RR\oplus w^{-1}$.

As shown in \cite{DonOrient}, the difference between the orientations
$o_{k}(\Omega,w')$ and $o_{k}(\Omega,w'')$ for $M_E^{\asd}$ is given by
\begin{equation}
\eps(w',w'')=(-1)^{(w'-w'')^2/4}.
\label{eq:CompareOrientation}
\end{equation}
where $w',w''\in H^2(X;\ZZ)$ are any two integral lifts of $w_2(\su(E))$.
\end{rmk}

\subsection{Geometric representatives and zero-section monopoles}
\label{subsec:CohomOnASDLink}
The stratum
$M^{\asd}_E\subset M_{W,E}$ of zero-section pairs is identified with
the moduli space of anti-self-dual connections on the
$\SO(3)$ bundle $\su(E)$.  Because the geometric representatives
$V(\beta)$ are pulled back by the map 
$\sC^*_{W,E}\to \sB^*_E$
given by $[A,\Phi]\mapsto [A]$,
the following computation of the intersection of the geometric
representatives with the stratum $M^{\asd}_E$ of
zero-section monopoles is clear:

\begin{lem}
\label{lem:ASDIntersection}
Let $E$ be a Hermitian two-plane bundle over a four-manifold $X$ with
$b^+(X)>0$ and generic Riemannian metric. Choose $c_1(E)\pmod{2}$ so that
$\su(E)$ does not admit a flat connection.
Let $z\in\AAA(X)$ have degree $2n_{p_1}$, where
$n_{p_1}\geq d_a$.  For a generic choice of geometric representatives, the
intersection of $\barV(z)$ with the strata of zero-section pairs in
$\barM_{W,E}$ is a finite number of generic points in $M^{\asd}_E$.

If $M^{\asd}_E$ is given its standard orientation then the number of points
in this intersection, counted with sign, is given by
$$
\#  (\barV(z)\cap \barM^{\asd}_E )
= 
\begin{cases}
D^{c_1(E)}_X(z)
& \text{if $n_{p_1}=d_a$},\\
0
&
\text{if $n_{p_1}>d_a$}.
\end{cases}
$$
\end{lem}

As we shall see in the following lemma,
it is important that the above intersection take place at generic points
in $M^{\asd}_E$.  
A neighborhood of a zero-section pair $[A,0]\in M_{W,E}$
can be described by the following Kuranishi model.

\begin{lem}
\label{lem:ASDKuranishi}
For any point $[A]\in M^{\asd}_E$,
there is a smoothly stratified diffeomorphism between
a neighborhood
of $[A,0]$ in $M_{W,E}$ and a neighborhood of zero
in $m^{-1}(0)/S^1_Z$, where $m$ is an $S^1_Z$-equivariant map
$$
m: T_A M^{\asd}_E\oplus \Ker D_{A,\vecvartheta}
\rightarrow
\Coker D_{A,\vecvartheta}.
$$
If $\Ind D_{A,\vecvartheta}>0$ then
for generic points $[A]\in M^{\asd}_E$, the cokernel of the Dirac 
operator vanishes for generic perturbations $\vecvartheta$.
\end{lem}

The cokernel of the perturbed Dirac operator $D_{A,\vecvartheta}$ vanishes
at generic points $[A]\in M^{\asd}_E$ because the map $A\mapsto
D_{A,\vecvartheta}$ from $\tM^{\asd}_E$ to the space Fredholm operators,
for a given index, is transverse to the jumping line strata.  As described
in \cite{Koschorke}, the `jumping line strata' are the strata of Fredholm
operators indexed by the dimension of their cokernels and the top stratum
consists of operators with vanishing cokernel.  Lemma
\ref{lem:ASDKuranishi} then describes the normal cone to $M^{\asd}_E$ at a
generic point $[A,0]$ as a cone on $\CC\PP^{n_a-1}$, where $\Ker
D_{A,\vecvartheta}\simeq \CC^{n_a}$.

We have described the geometric representative $\barV(\beta)$ near
the anti-self-dual moduli space; $\barW(x)$ can be described as follows.

\begin{lem}
\label{lem:CohomOnASDNormalSlice}
When restricted to the link in the 
normal cone of $M^{\asd}_E$ in $M_{W,E}$
at a generic point $[A,0]\in M^{\asd}_E$, the geometric representative
$\barW(x)$ is Poincare dual to $2h$, where 
$h\in H^2( (\Ker D_{A,\vecvartheta}\backslash \{0\})/S^1_Z;\ZZ)$ 
is the positive generator.
\end{lem}

\begin{rmk}
Lemma \ref{lem:CohomOnASDNormalSlice} shows that $W(x)$
will have non-trivial intersection with the normal cone
of any generic point in $M^{\asd}_E$.  Thus, the closure
of $W(x)$ in $M_{W,E}$
will contain all generic points and thus all points in
$M^{\asd}_E$.
\end{rmk}

Let $\barM^{\asd}_E$ denote the closure of $M^{\asd}_E$ in
$\barM_{W,E}$; note that this may properly contain the closure
$M^{\barasd}_E$ of $M^{\asd}_E$ in $IM^{\asd}_E$.

\begin{defn} 
The {\em link of $\barM^{\asd}_E$\/} in $\barM_{W,E}$ is given by
$$
\bL^{\asd}_\varepsilon
:=
\{[A,\Phi]\in \barM_{W,E}: \|\Phi\|_{L^2}^2=\varepsilon^2\}.
$$
\end{defn}

It is a simple matter to show that the map $\Phi\mapsto \|\Phi\|_{L^2}^2$
is continuous on $\barM_{W,E}$ and smooth on each stratum.  Thus, for
generic values of $\varepsilon>0$, the link $\bL^{\asd}_\varepsilon$ is a
smoothly stratified, codimension-one subspace of $\barM_{W,E}$.  The
intersection of $\bL^{\asd}_\varepsilon$ with an approriate number of
generic geometric representatives is then a finite number of points
which can be calculated using Lemmas 
\ref{lem:ASDIntersection} and \ref{lem:CohomOnASDNormalSlice}.

\begin{lem}
\label{lem:LinkOfASD}
\cite{FL2}
Let $E$ be a Hermitian two-plane bundle over a four-manifold $X$ with
$b^+(X)>0$ and generic Riemannian metric. Choose $c_1(E)\pmod{2}$ so that
$\su(E)$ does not admit a flat connection. Let $n_{p_1}$ and $n_{c_1}$ be
non-negative integers such that $n_{p_1}+n_{c_1}=d_a+n_a-1$.  Suppose
$z\in\AAA(X)$ has degree $2n_{p_1}\geq 2d_a$.  If $M_{W,E}$ is given the
orientation $O_{k}^{\asd}(\Omega,c_1(E))$, where $k = -\quarter
p_1(\su(E))$, then there is a positive constant $\varepsilon_0$ such that
for generic $\varepsilon<\varepsilon_0$ we have
$$
\# \left(\barV(z)\cap \barW(x^{n_{c_1}})\cap \bL^{\asd}_\varepsilon \right) 
 = 
\begin{cases}
2^{n_a-1}D^{c_1(E)}_X(z)
& \text{if $n_{p_1}=d_a$},\\
0
&
\text{if $n_{p_1}>d_a$}.
\end{cases}
$$
\end{lem}

\subsection{Links of the strata of reducible monopoles}
\label{subsec:LinkOfReducible}
To describe the geometric representatives in a neighborhood of the
reducible monopoles, $M^{\red}_{W,E,L_1}$, it does not suffice to produce a
Kuranishi model at a generic point.  Neither of the geometric
representatives, $V(\beta)$, $W(x)$ intersects $M^{\red}_{W,E,L_1}$ in a
set of the expected codimension so we cannot use them to cut down to a set
of generic points as we did with the stratum of zero-section monopoles.
Instead, we must give a global description of the link of
$M^{\red}_{W,E,L_1}$ in $M_{W,E}^{*,0}$. We may assume without loss of
generality that $M^{\red}_{W,E,L_1}$ contains no zero-section solutions.

Even in the case where $M^{\red}_{W,E,L_1}$ is in the top level $M_{W,E}$,
the problem of defining a link is non-trivial when the dimension of
$M^{\red}_{W,E,L_1}$ is positive.  The techniques we employ in \cite{FL2}
follow the ideas of Atiyah and Singer for stabilizing index bundles
\cite{ASFamily,DK}. Related methods have also been used in a variety of
recent applications of Gromov and Seiberg-Witten invariants (including
those of \cite{Brussee,Furuta,LiLiu,LiTian,RuanSW,RuanGW}, for example)
which essentially involve `excess intersection theory' in situations where
transversality cannot be achieved by `generic parameter' arguments via the
Sard-Smale theorem.

In this subsection, we sketch our construction of the link of
$M^{\red}_{W,E,L_1}$ in $M_{W,E}^{*,0}$ when these reducibles lie in the top
level \cite{FL2}.  Let $(A,\Phi)$ represent a point in
$M^{\red}_{W,E,L_1}$ and recall that
$\sD_{A,\Phi}=d^{0,*}_{A,\Phi}+d^1_{A,\Phi}$. Let 
$$
\bE 
:=
L^2_{k-1}(\Lambda^+\otimes\su(E))\oplus L^2_{k-1}(W^-\otimes E)
$$ 
and let $\fS:\tsC_{W,E}\to\bE$ be the $\ssG_E$-equivariant map 
defined by the $\PU(2)$ monopole equations \eqref{eq:PT}, so
$d^1_{A,\Phi}=(D\fS)_{A,\Phi}$. It is convenient to temporarily pass to an
$S^1$-equivariant setting, so let
$$
\ssC_{W,E} := \tsC_{W,E}/\sG_E,
$$
and note that $\sC_{W,E} = \ssC_{W,E}/S_Z^1 = \ssC_{W,E}/S_{L_2}^1$. We
then have 
$$
\cM_{W,E} := \fS^{-1}(0)\cap\ssC_{W,E},
$$ 
with quotient $M_{W,E} = \cM_{W,E}/S_Z^1 = \cM_{W,E}/S_{L_2}^1$.  If
$[A,\Phi]$ is a point in $M^{\red}_{W,E,L_1}$, the stabilizer
$\Stab_{A.\Phi}$ of the pair $(A,\Phi)$ is $S_{L_2}^1$ in $\ssG_E$ but is
trivial in $\sG_E$.

If $(A,\Phi)$ represents a point $[A,\Phi]\in M^{\red}_{W,E,L_1}$,
the full elliptic deformation complex $d_{A,\Phi}^\bullet$ of
\eqref{eq:DeformationComplex} for the $\PU(2)$ monopole equations splits
into {\em tangential deformation complex\/}, $d_{A,\Phi}^{\bullet,t}$, and
{\em normal deformation complex\/}, $d_{A,\Phi}^{\bullet,n}$ (see
\cite{FL2}). The tangential deformation complex is isomorphic to the
elliptic deformation complex for the $\U(1)$ monopole equations
\eqref{eq:U1Monopole}. The rolled-up elliptic
deformation complex $\sD_{A,\Phi}=d^{0,*}_{A,\Phi}\oplus d^1_{A,\Phi}$ also
splits, of course, into tangential and normal rolled-up deformation
complexes: $\sD_{A,\Phi}=\sD^t_{A,\Phi}\oplus\sD^n_{A,\Phi}$, with
$\Coker\sD^t_{A,\Phi} = 0$ and
$$
\Ker\sD^t_{A,\Phi} \simeq T_{A,\Phi}M^{\red}_{W,E,L_1}
\quad\text{and}\quad
\Ker\sD^n_{A,\Phi} \simeq \Ker\sD_{A,\Phi}/T_{A,\Phi}M^{\red}_{W,E,L_1}.
$$
Let $\Pi_{A,\Phi}$ denote the $L^2$ orthogonal projections onto the
subspaces
$$
\Coker d_{A,\Phi}^1 
\simeq 
\Coker\sD_{A,\Phi}
\simeq 
\Coker\sD_{A,\Phi}^n
\simeq 
\Coker d_{A,\Phi}^{1,n}, 
$$
noting that $\Coker d_{A,\Phi}^{0,*}=\Ker d_{A,\Phi}^0=0$. The
Kuranishi model of a neighborhood in $M_{W,E}$ of a point $[A,\Phi]\in
M^{\red}_{W,E,L_1}$ is given by
\begin{gather}
\bga: \sO_{A,\Phi}\subset \Ker\sD_{A,\Phi} \to \ssC_{W,E}, 
\label{eq:KuranishiModel}
\\
\bvarphi:\sO_{A,\Phi}\subset \Ker\sD_{A,\Phi}\to \Coker\sD_{A,\Phi}, 
\notag
\end{gather}
where $\sO_{A,\Phi}$ is an $S^1_{L_2}$ invariant open neighborhood of the
origin in $\Ker\sD_{A,\Phi}=\Ker\sD_{A,\Phi}^t\oplus\Ker\sD_{A,\Phi}^n$, 
$\bga$ is an $S^1_{L_2}$-equivariant embedding, and $\bvarphi$ is a
smooth $S^1_{L_2}$-equivariant map. The map $\bga$ descends to a
smoothly stratified diffeomorphism from $\bvarphi^{-1}(0)/S^1_{L_2}$ onto   
an open neighborhood of $[A,\Phi]$ in $M_{W,E}$. The
obstruction map $\bvarphi$ is given by $\Pi\circ\fS\circ\bga$.

Since the construction of the link of $M^{\red}_{W,E,L_1}$ in
$M_{W,E}^{*,0}$ is complicated in general, it is helpful to begin by
considering some simple special cases.  When $M^{\red}_{W,E,L_1}$ is
zero-dimensional, links in $M_{W,E}^{*,0}$ of the points of
$M^{\red}_{W,E,L_1}$ are defined by the Kuranishi model
\eqref{eq:KuranishiModel}: The link of a point $[A,\Phi]$ is simply given
by the $S^1_{L_2}$ quotient of the zero-locus of $\bvarphi$ in an
$\eps$-sphere around the origin in 
$\Ker\sD_{A,\Phi}$.

For the remainder of this subsection we assume that $M^{\red}_{W,E,L_1}$
may be positive-dimensional. If $\Coker\sD$ vanishes along
$M^{\red}_{W,E,L_1}$, then $\Ker\sD^n$ is a finite-rank, 
$S^1_{L_2}$-equivariant vector bundle over $M^{\red}_{W,E,L_1}$ with fibers
$\Ker\sD^n_{A,\Phi}$ over points $[A,\Phi]\in M^{\red}_{W,E,L_1}$.  There
is an $S^1_{L_2}$-equivariant diffeomorphism $\bvarphi$ from an open
neighborhood $\sO$ of the zero section $M^{\red}_{W,E,L_1}\subset\Ker\sD^n$
and an open neighborhood of $M^{\red}_{W,E,L_1}$ in $\cM_{W,E}$.

More generally, if the cokernel of $\sD_{A,\Phi}$ has constant rank as
$[A,\Phi]$ varies in $M^{\red}_{W,E,L_1}$ (that is, no spectral flow
occurs), then $\Ker\sD^n$ and $\Coker\sD$ both define finite-rank,
$S^1_{L_2}$-equivariant vector bundles over $M^{\red}_{W,E,L_1}$:
\begin{equation}
\begin{CD}
\Ker\sD^n @.                     @.                @.@. \Coker\sD     \\
          @.   {\pi_k}\searrow   @.                @. \swarrow{\pi_c} \\
          @.                     @. \quad M^{\red}_{W,E,L_1}
\end{CD}
\label{eq:KernelCokernelBundle}
\end{equation}
Let $2\nu$ be the least positive eigenvalue of the Laplacian $\De_{A,\Phi}
:= \sD_{A,\Phi}\sD_{A,\Phi}^*$ as $[A,\Phi]$ varies along the compact manifold
$M^{\red}_{W,E,L_1}$ and let $\Pi_{\nu;A,\Phi}$ denote the $L^2$ orthogonal
projection from $\bE$ onto the subspace spanned by the eigenvectors
of $\De_{A,\Phi}$ with eigenvalue less than $\nu$.  
The vector bundle $\Coker\sD$ over $M^{\red}_{W,E,L_1}$ then extends to a
vector bundle $\Xi_\nu:= \Ker(\id-\Pi_\nu)\circ\De =
\Coker(\id-\Pi_\nu)\circ\sD$ of the same rank over 
an open neighborhood of $M^{\red}_{W,E,L_1}$ in $\csC_{W,E}$.
The {\em obstruction section\/} $\bvarphi$ over $\sO\subset \Ker\sD^n$ of
the vector bundle
\begin{equation}
\bga^*\Xi_\nu \to \Ker\sD^n
\label{eq:ObstructionProjection}
\end{equation}
is given by $\bvarphi := \Pi_\nu\circ\fS\circ\bga$ on
$\sO\subset\Ker\sD^n$, where the $S^1_{L_2}$-equivariant embedding 
$\bga:\sO\to\csC_{W,E}$ gives a 
diffeomorphism from an open neighborhood $\sO$ of the zero section
$M^{\red}_{W,E,L_1}$ in $\Ker\sD^n$ onto an open neighborhood of
$M^{\red}_{W,E,L_1}$ in the $S^1_{L_2}$ invariant {\em thickened moduli
space\/}
$$
\cM_{W,E,L_1}(\Xi_\nu) := ((\id-\Pi_\nu)\circ\fS)^{-1}(0) \subset\csC_{W,E}.
$$
Then $\bga$ descends to a smoothly stratified diffeomorphism from the zero
locus
\begin{equation}
\bvarphi^{-1}(0)/S^1_{L_2}\subset\Ker\sD^n/S^1_{L_2}
\label{eq:ZeroLocus}
\end{equation}
containing $M^{\red}_{W,E,L_1}$ onto an open neighborhood of
$M^{\red}_{W,E,L_1}$ in $M_{W,E}$.  On the complement of the zero section
$M^{\red}_{W,E,L_1}\subset
\Ker\sD^n$, the $S^1_{L_2}$ quotient of the projection
\eqref{eq:ObstructionProjection} given by
\begin{equation}
\bga^*\Xi_\nu/S^1_{L_2} \to \Ker\sD^n/S^1_{L_2},
\label{eq:ObstructionBundle}
\end{equation}
is a vector bundle. The homology class of the zero locus 
\eqref{eq:ZeroLocus} of the obstruction map can be calculated from the
Euler class of the vector bundle \eqref{eq:ObstructionBundle} or,
equivalently, from that of
$$
\pi_k^*\Coker\sD/S^1_{L_2} \to \Ker\sD^n/S^1_{L_2},
$$
as is easily seen.

In general, though, one cannot guarantee that $\Coker\sD$ will
either vanish or have constant rank.  Let
$\tM^{\red}_{W,E,L_1}\subset\tsC_{W,E}$ be the pre-image of
$M^{\red}_{W,E,L_1}$ under the projection from the pre-configuration space
$\tsC_{W,E}$ onto the quotient $\csC_{W,E} = \tsC_{W,E}/\sG_E$.  Because
$M^{\red}_{W,E,L_1}$ is compact, we can construct a finite family of 
gauge-equivariant `stabilizing maps' from $\tM^{\red}_{W,E,L_1}$ to
$\bE$ such that 
\begin{itemize}
\item
The image $\Xi_{A,\Phi}$
of these maps at $(A,\Phi)\in\tM^{\red}_{W,E,L_1}$
spans $\Coker \sD_{A,\Phi}$,
\item
The subspace $\Xi_{A,\Phi}\subset \bE$ is $S^1_{L_2}$ invariant,
\item
The dimension of $\Xi_{A,\Phi}$ is constant for all pairs
$(A,\Phi)\in \tM^{\red}_{W,E,L_1}$.
\end{itemize}
The subspaces $\Xi_{A,\Phi}$ then fit together to form an 
$S^1_{L_2}$-equivariant vector bundle $\Xi$ over $M^{\red}_{W,E,L_1}$,
which extends to 
an $S^1_{L_2}$-equivariant vector bundle $\Xi$ over an open neighborhood of
$M^{\red}_{W,E,L_1}$ in $\csC_{W,E}$. Let $\Pi_{\Xi;A,\Phi}$ denote the
$L^2$ orthogonal projection from $\bE$ onto the subspace $\Xi_{A,\Phi}$.  The
properties of the stabilizing sections ensure that the space
$$
N_{W,E,L_1}(\Xi) := \Ker(\id-\Pi_\Xi)\circ\sD^n
$$ 
is a vector bundle over
$M^{\red}_{W,E,L_1}$ with fibers which are closed under the $S^1_{L_2}$
action:
\begin{equation}
\begin{CD}
N_{W,E,L_1}(\Xi) @.              @.                @.@. \quad\Xi        \\
          @.   {\pi_N}\searrow   @.                @. \swarrow{\pi_\Xi} \\
          @.                     @. \quad M^{\red}_{W,E,L_1}
\end{CD}
\label{eq:NormalObstructionBundle}
\end{equation}
The bundle $\Xi$ plays the role of $\Xi_\nu$ while
$N_{W,E,L_1}(\Xi)$ plays that of $\Ker\sD^n$ in the simpler case 
\eqref{eq:KernelCokernelBundle} where the
cokernel of $\sD$ has constant rank along $M^{\red}_{W,E,L_1}$.  In
\cite{FL2} we construct a smooth,  
$S^1_{L_2}$ invariant thickened moduli space,
$$
\cM_{W,E,L_1}(\Xi) := ((\id-\Pi_\Xi)\circ\fS)^{-1}(0) \subset\csC_{W,E},
$$
using the stabilizing bundle $\Xi$. Then $N_{W,E,L_1}(\Xi)$ is the
$S^1_{L_2}$-equivariant normal bundle of the smooth submanifold
$M^{\red}_{W,E,L_1}\subset \cM_{W,E,L_1}(\Xi)$, recalling that
$M^{\red}_{W,E,L_1}$ is the fixed-point set of $S^1_{L_2}$. 

The equivariant tubular neighborhood theorem provides 
an $S^1_{L_2}$-equivariant diffeomorphism 
$\bga:\sO\to\csC_{W,E}$ from
an open neighborhood $\sO$ 
of the zero-section $M^{\red}_{W,E,L_1}\subset N_{W,E,L_1}(\Xi)$
onto an open neighborhood of the submanifold
$M^{\red}_{W,E,L_1}\subset\cM_{W,E,L_1}(\Xi)$
which covers the identity on $M^{\red}_{W,E,L_1}$.
The map $\bga$ then descends to a smoothly stratified diffeomorphism
from the zero locus $\bvarphi^{-1}(0)/S^1_{L_2}$ in
$N_{W,E,L_1}(\Xi)/S^1_{L_2}$ onto an open neighborhood of
$M^{\red}_{W,E,L_1}$ in the actual moduli space, $M_{W,E}$, where 
$$
\bvarphi := \Pi_\Xi\circ\fS\circ\bga
$$ 
is a section over $\sO\subset N_{W,E,L_1}(\Xi)$ of the 
$S^1_{L_2}$-equivariant vector bundle
$$
\bga^*\Xi \to N_{W,E,L_1}(\Xi).
$$
As in the constant rank case, this descends to a vector bundle 
$$
\bga^*\Xi/S^1_{L_2} \to N_{W,E,L_1}(\Xi)/S^1_{L_2}
$$
on the complement of the
zero section, $M^{\red}_{W,E,L_1}\subset N_{W,E,L_1}(\Xi)/S^1_{L_2}$, whose
Euler class may be computed from
$$
\pi_N^*\Xi/S^1_{L_2} \to N_{W,E,L_1}(\Xi)/S^1_{L_2}.
$$
While the bundle $\bga^*\Xi$ given by this
restriction to the complement of the zero section
is trivial --- because it is spanned by the stabilizing sections ---
the quotient $\bga^*\Xi/S^1_{L_2}$ has a non-trivial Euler class.

\begin{defn}
Let $N_{W,E,L_1}^\varepsilon(\Xi)$ denote the sphere bundle
of fiber vectors of length $\varepsilon$ and 
let $\PP N_{W,E,L_1}(\Xi)=N_{W,E,L_1}^\varepsilon(\Xi)/S^1_{L_2}$.
The {\em link of the stratum $M^{\red}_{W,E,L_1}\subset M_{W,E}$\/}
of reducible pairs is given by 
$$
\bL_{W,E,L_1}
:=
\bga\left(\bvarphi^{-1}(0)\cap N_{W,E,L_1}^\varepsilon(\Xi)\right)/S^1_{L_2}
$$
and thus
$$
[\bL_{W,E,L_1}] 
= 
e\left(\bga^*\Xi/S^1_{L_2}\right)\cap [\PP N_{W,E,L_1}(\Xi)]
$$
is its homology class.
\end{defn}

\begin{rmk}
The orientation given to $\bL_{W,E,L_1}$ by the orientation
on $M^{\red}_{W,E,L_1}$ from the homology orientation $\Omega$
and the complex structure on the fibers of $N_{W,E,L_1}(\Xi)$ (from
the $S^1_{L_2}$ action) is equivalent to the orientation given by
$O^{\asd}_{k}(\Omega,L_2\otimes L_1^*)$ (see \cite{FL2}).
\end{rmk}

\subsection{Reduction formulas for Donaldson invariants:
$\U(1)$ monopoles in the top Uhlenbeck level} In this subsection we
describe some of our results from \cite{FL2}, where we compute Donaldson
invariants in terms of Seiberg-Witten invariants when the $\U(1)$ monopoles
in $\barM_{W,E}$ lie only in the top level $M_{W,E}$.

\begin{defn}
\label{defn:SetOfReducibles}
The set of moduli spaces of $\U(1)$ monopoles 
contained in the top level $M_{W,E}$ is enumerated by
$$
\Red(W,E)
:=
\{L_1\in H^2(X;\ZZ): M^{\red}_{W,E,L_1}\neq \emptyset
\text{ and } (2L_1 - c_1(E))^2 = p_1(\su(E))\}. 
$$
The set of moduli spaces of $\U(1)$ monopoles 
contained in the compact space of ideal $\PU(2)$ monopoles $IM_{W,E}$ is
enumerated by 
\begin{multline*}
\overline\Red(W,E)
:=
\{L_1\in H^2(X;\ZZ): M^{\red}_{W,E_{-\ell},L_1}\neq \emptyset \text{ and } 
\\
(2L_1 + c_1(E))^2 = p_1(\su(E_{-\ell}))+4\ell, \ell\in\ZZ_{\ge 0}\}, 
\end{multline*}
where $c_1(E_{-\ell}) = c_1(E)$ and $c_2(E_{-\ell}) = c_2(E)-\ell$.
\end{defn}

Note that $2L_1 - c_1(E) = K-F$, where $K=c_1(W^+\otimes L_1)$ and
$F=c_1(W^+)+c_1(E)$.  The compactification $\barM_{W,E}$ may be a proper
subset of $IM_{W,E}$.  If the reducibles in $\barM_{W,E}$ appear only in
the top Uhlenbeck level $M_{W,E}$ then $\barM^{*,0}_{W,E}$ serves as a
cobordism between the link $\bL^{\asd}_\eps$ of the anti-self-dual moduli
space $M_E^{\asd}$ and the links $\bL_{W,E,L_1}$ of the strata of
reducibles $M_{W,E,L_1}^{\red}$.  This gives the following formula:

\begin{thm}
\label{thm:CompactReductionFormula}
\cite{FL2}
Let $E$ be a Hermitian two-plane bundle over a four-manifold $X$ with
$b^+(X)>0$ and generic Riemannian metric. Choose $c_1(E)\pmod{2}$ so that
$\su(E)$ does not admit a flat connection. Suppose $z\in\AAA(X)$ has
degree $2d_a$.  If $\Red(W,E)=\overline\Red(W,E)$, so the reducible
$\PU(2)$ monopoles in $\barM_{W,E}$ appear only in the highest Uhlenbeck
level, then
\begin{equation}
 2^{n_a-1}  D^{c_1(E)}_X(z) 
=
-\sum_{L_1\in \Red(W,E)} (-1)^{L_1^2}
\langle \mu_{p_1}(z)\smile\mu_{c_1}(x^{n_a-1}),[\bL_{W,E,L_1}]\rangle.
\label{eq:CompactReductionFormula1}
\end{equation}
\end{thm}

The sign $(-1)^{L_1^2}$ in \eqref{eq:CompactReductionFormula1} comes from
the parity change $\eps(c_1(E),L_2\otimes L_1^*)$ of
\eqref{eq:CompareOrientation}, noting that $c_1(E) = L_1+L_2$.

The restriction of the cohomology classes $\mu_{p_1}(\beta)$ and $\mu_{c_1}(x)$
to $\bL_{W,E,L_1}$ are computed in \cite{FL2} in terms of
the hyperplane class on $M^{\red}_{W,E,L_1}$ and the generator of
the cohomology of the fiber of $\PP N_{W,E,L_1}(\Xi)$.  
The Euler class, $e(\bga^*\Xi/S^1_{L_2})$, can also be expressed in these
terms. {}From the Atiyah-Singer index theorem for families, one
can compute the Segre classes of the bundle $N_{W,E,L_1}(\Xi)$
under the assumption $b^1(X)\leq 1$.  If $b^1(X)>1$ the computation
is still possible in principle, but becomes unmanageable in practice.
To describe the results of these computations,
we introduce some standard expressions to describe certain
constants arising in our reduction formula:

\begin{defn}
\cite[\S 8.96]{GR}
\label{defn:JacobiPoly}
The {\em Jacobi polynomials\/} are defined by
$$
P^{(a,b)}_n(x) 
:= 
2^{-n}
\sum_{m=0}^n {\binom{n+a}{m}}{\binom{n+b}{n-m}}
(x-1)^{n-m}(x+1)^{m}.
$$
\end{defn}

Functional relations, relations with other special functions, and the
generating function for the Jacobi polynomials can be found in
\cite[pp. 1034--1035]{GR}. Recall that $\fs_0 = (\rho,W)$ is a
choice of fixed \spinc structure on $X$. For line bundles $L_1
\in H^2(X;\ZZ)$, we denote $\fs_0\otimes L_1 := (\rho,W\otimes L_1)$.

\begin{thm}
\label{thm:DegreeZeroFormula}
\cite{FL2}
Let $E$ be a Hermitian two-plane bundle over a four-manifold $X$ with
$b^+(X)>0$, $b^1(X)\leq 1$, and generic Riemannian metric. Choose
$c_1(E)\pmod{2}$ so that $\su(E)$ does not admit a flat connection. Let
$n_{p_1}+n_{c_1}=d_a+n_a-1$, where $n_{p_1}$, $n_{c_1}$ are non-negative
integers.  For the stratum of reducible solutions $M^{\red}_{W,E,L_1}$
contained in the highest level of $\barM_{W,E}$, a generator $x\in
H_0(X;\ZZ)$, classes $\beta_1,\dots,\beta_{n_{p_1}}\in H_2(X;\QQ)$, and
integers $0\le m\le
\half n_{p_1}$, we have
\begin{align}
\label{eq:DegreeZeroPairing}
&\left\langle
\mu_{p_1}(\beta_1\cdots\beta_{n_{p_1}-2m}x^m)\smile\mu_{c_1}(x^{n_{c_1}}),
[\bL_{W,E,L_1}]
\right\rangle\\
\notag
&\qquad =
(-1)^m
2^{-n_{p_1}+d_s}
C_{W,E,L_1}(n_{p_1},n_{c_1})
\SW(\fs_0\otimes L_1)
\prod_{i=0}^{n_{p_1}-2m}\langle 2L_1 - c_1(E),\beta_i\rangle
\end{align}
where, for $I=n_{p_1}-n^\Lambda_s-d_s$ and $J=n_{c_1}-d_s$, the constants 
$n^\Lambda_s$ and $C_{K,F}$ are given by  
\begin{align*}
n^\Lambda_s(\su(E)) &:= -p_1(\su(E)) -\half (e(X)+\sigma(X)),
\\
C_{W,E,L_1}(n_{p_1},n_{c_1})
&:=  
P^{(I,J)}_{d_s}(0)
=
2^{-d_s}\sum_{u=0}^{d_s}
(-1)^u{\binom {n_{c_1}} u}{\binom {n_{p_1}-n^\Lambda_s}{d_s-u}}.
\end{align*}
\end{thm}

\begin{rmk}
\begin{enumerate}
\item Note that $2L_1 - c_1(E) = K-F$, where $K=c_1(W^+\otimes L_1)$
and $F=c_1(W^+)+c_1(E)$, and that the polynomial $C_{W,E,L_1}(\cdot)$ only
depends on the classes $K$ and $F$ (together with the Euler characteristic
and signature of $X$).
\item The constant $n^\Lambda_s$ is the index of the elliptic complex on
$\Omega^\bullet(L_1\otimes L_2^*)$ induced by homotoping the normal
deformation complex at a reducible pair, determined by the reduction $E =
L_1\oplus L_2$, to a diagonal complex. 
\item If $d_s=0$ we have $P^{(I,J)}_0(0)=1$ and so for manifolds of
Seiberg-Witten simple type, the constant $C_{W,E,L_1}(n_{p_1},n_{c_1})$ is
not interesting.  It should, however, prove vital in understanding the
relation between the Donaldson and Seiberg-Witten invariants for any manifolds
which are not of simple type.
\end{enumerate}
\end{rmk}

Combining Theorems \ref{thm:CompactReductionFormula} and 
\ref{thm:DegreeZeroFormula} yields:

\begin{cor}
\label{cor:DegreeZeroFormula}
\cite{FL2}
Let $E$ be a Hermitian two-plane bundle over a four-manifold $X$ with
$b^+(X)>0$, $b^1(X)\leq 1$, and generic Riemannian metric. Choose
$c_1(E)\pmod{2}$ so that $\su(E)$ does not admit a flat connection.  Let
$x\in H_0(X;\ZZ)$ be a generator, let $\beta_1,\dots,\beta_{d_a}\in
H_2(X;\QQ)$, and suppose
$$
z = \beta_1\cdots\beta_{d_a-2m}x^m \in \AAA(X),
$$
for $0\le m \le \half d_a$. If $\Red(W,E)=\overline\Red(W,E)$, so
reducible $\PU(2)$ monopoles in $\barM_{W,E}$ appear only in the highest
level $M_{W,E}$, then the following holds:
\begin{multline*}
-2^{n_a-1}D^{c_1(E)}_X(z) 
=
\sum_{L_1\in \Red(W,E)}
(-1)^{L_1^2}(-1)^m
2^{-d_a+d_s(c_1(W^+\otimes L_1))} \\
\times
C_{W,E,L_1}(d_a,n_a-1)\SW(\fs_0\otimes L_1)
\prod_{i=0}^{d_a-2m}\langle 2L_1 - c_1(E),\beta_i\rangle,
\end{multline*}
where $C_{W,E,L_1}(d_a,n_a-1)$ is defined in
Theorem \ref{thm:DegreeZeroFormula}. If $X$ has Seiberg-Witten
simple type then
$$
D^{c_1(E)}_X(z) 
=
\sum_{L_1\in \Red(W,E)}
(-1)^{L_1^2}(-1)^{m-1}2^{1-d_a-n_a}
\SW(\fs_0\otimes L_1)
\prod_{i=0}^{d_a-2m}\langle 2L_1 - c_1(E),\beta_i\rangle.
$$
\end{cor}

The formula in Corollary \ref{cor:DegreeZeroFormula} differs what one might
expect from equations \eqref{eq:KMFormula} and \eqref{eq:WittenFormula} as it
contains terms of the form
$$
\langle 2L_1-c_1(E),\beta_i\rangle
=
\langle K-F,\beta_i\rangle,
$$
where $K=c_1(W^+\otimes L_1)$ and $F=c_1(E)+c_1(W^+)$, rather than the
terms $\langle K,\beta_i\rangle$.  In addition, the power $L_1^2$ of $-1$
does not match the exponent $\half(w^2+wK)$ given in \eqref{eq:KMFormula}
for any obvious choice of line bundle $w$ over $X$.

As shown by our examples in \cite{FL2}, the condition
$\Red(W,E)=\overline\Red(W,E)$ puts severe restrictions
on the class $F$ and the intersections $FK_r$, where
the $K_r$ are basic classes.  Under these restrictions,
combinatorial identities give a cancellation of the factors
of $F$ in the formula of Corollary \ref{cor:DegreeZeroFormula}.
One sees from these examples that one should not assume that the terms
$$
(-1)^{\half(w^2+wK_r)}\exp(Q/2)\SW(K_r)e^{K_r}
$$
in \eqref{eq:WittenFormula} translate directly into
values for pairings with the link of the reducible $M^{\red}_{W,E,L_1}$
when $K=c_1(W^+\otimes L_1)$.  In the sum over all links, there can be many
cancellations between terms contributed by different links. We illustrate
the use of Corollary \ref{cor:DegreeZeroFormula} below;
see \cite{FL2} for further examples.

\begin{exmp}
\label{exmp:Elliptics}
\cite{FL2}
We use Corollary \ref{cor:DegreeZeroFormula} to
calculate the first non-trivial Donaldson polynomial of the elliptic
surface $E(n)$ with Euler characteristic $e(E(n))=12n$ and signature
$\sigma(E(n))=-8n$.  Let $f\in H^2(E(n);\ZZ)$ denote the fiber class of the
elliptic fibration.  For suitable perturbations, the only non-empty
Seiberg-Witten moduli spaces correspond to \spinc structures with
$$
K_r := c_1(W^+\otimes L_{1,r})=(n-2-2r)f, \qquad r=0,\dots,n-2.
$$
The Seiberg-Witten invariants of the \spinc structures with these
classes are given by (see, for example, \cite{FSRationalBlowDown}):
$$
\SW(K_r)=(-1)^r{\binom {n-2} r}, \qquad r=0,\dots,n-2.
$$
Because $p_1(\su(E)) = (L_1-L_2)^2 = (K_r-F)^2$, where 
$E=L_{1,r}\oplus (\det E)\otimes L_{1,r}^*$,
we can ensure that all the reducibles are in the same
level (and make this the top level) by requiring that $K_rF=0$.  Then
$p_1(\su(E))=(K_r-F)^2=F^2$.  Since $(1+b^+(E(n)))=2n$, we find that
\begin{align*}
d_a(\su(E)  &= -F^2-\threehalf(2n)=-F^2-3n, \\
n_a(\su(E)) &=\quarter(2F^2+8n)=\half F^2+2n.
\end{align*}
Thus, to obtain $d_a\geq 0$ and $n_a>0$, we impose the constraint $-4n<
F^2\leq -3n$. Note that as $K_r$ is characteristic and $K_rF=0$, we must
have $F^2$ even. Applying Corollary \ref{cor:DegreeZeroFormula} with $\beta
\in H_2(X;\ZZ)$ we find, after some calculation, that
$$
D^F_X(\beta^{n-2j-2m}x^m)=
\begin{cases}
0 & \text{if $j>1$ or $m>0$}, \\
-(n-2)!\langle f,\beta\rangle^{n-2} & \text{if $j=m=0$},
\end{cases}
$$
in agreement with the results of \cite{FSRationalBlowDown,KMStructure}.
\end{exmp}


\section{Gluing $\PU(2)$ monopoles and the $\PU(2)$ monopole analogue
of the Kotschick-Morgan conjecture} 
\label{sec:GluingKMConjecture}
The problems involved in computing intersection
numbers for the link of a family of lower-level reducibles are similar
to those encountered in attempts to prove the Kotschick-Morgan conjecture
\cite{KoM}. In this section we first discuss the Kotschick-Morgan
conjecture for Donaldson invariants, describe its analogue for
pairings with links of lower-level moduli spaces 
of $\U(1)$ monopoles in the
Uhlenbeck compactification of the moduli space of $\PU(2)$ monopoles, and
outline how a resolution of this analogue should lead in turn to a proof of
Witten's conjecture.

\subsection{The Kotschick-Morgan conjecture for Donaldson invariants}
\label{subsec:KMConjecture}
The conjecture of Kotschick and Morgan for Donaldson invariants of
four-manifolds $X$ with $b^+(X)=1$ gives a prediction of how the Donaldson
invariants vary when the underlying Riemannian metric changes. More
precisely, it asserts that the invariants computed using metrics lying in
different chambers of the positive cone of $H^2(X;\RR)/\RR^*$ differ by
terms depending only the homotopy type of $X$ \cite{KoM}. The definition of
the Donaldson invariants requires a choice of Riemannian metric on $X$ and
they are only independent of this choice when $b^+(X)>1$.

The Donaldson invariants of a manifold with $b^+(X)=1$ are not independent
of the metric because the cobordism formed by taking the moduli space of
connections anti-self-dual with respect to elements of a path of metrics
may contain reducible anti-self-dual connections.  The Donaldson cohomology
classes evaluate non-trivially on the links of these reducible connections,
so the values of the Donaldson polynomial given by the metrics at the ends
of this path will differ by the pairing of the top power of the cohomology
classes with these links.  Directly evaluating such pairings or even
showing that they depend only on homotopy data is a difficult problem when
the reducible connection lies in a lower level of the Uhlenbeck
compactification.  The conjecture of \cite{KoM} asserts that these pairings
only depend on homotopy data: this has been verified for reducibles in the
strata $M^{\asd}_{E_{-\ell}}(X)\times \Sym^\ell(X)$ when $\ell\le 2$
\cite{DonHCobord,Kotschick,KoM,Leness,Y} and for much higher $\ell$ when $X$
is algebraic \cite{EG,FrQ}.

Motivated by related work of L. G\"ottsche on the Kotschick-Morgan
conjecture for Donaldson invariants \cite{Goettsche} and by Fintushel and
Stern on the general blowup formula \cite{FSBlowupD}, Pidstrigach and
Tyurin suggested that the conjecture of Witten should then follow by
calculations --- analogous to those of G\"ottsche --- from the
Kotschick-Morgan conjecture for $\PU(2)$ monopoles \cite{PTLocal}. In the
case of $\PU(2)$ monopoles there are further complications, not present in
Donaldson theory, due in part to the many additional obstructions to gluing
$\PU(2)$ monopoles.

\subsection{$\PU(2)$ monopoles: Gluing and ungluing}
\label{subsec:PTGluing}
The cobordism scheme requires the use of analogues of
Taubes' gluing maps to parametrize neighborhoods of moduli
spaces of $\U(1)$ monopoles
lying at the Uhlenbeck boundary of the moduli space of $\PU(2)$
monopoles and in particular, to construct links of these singularities.

In our articles \cite{FL3,FL4} we first construct approximate gluing maps
--- giving approximate solutions to the $\PU(2)$ monopole equations ---
by grafting anti-self-dual connections from the four-sphere, which are
concentrated at the north pole, onto a
background $\PU(2)$ monopole at distinct points which are allowed to
vary. We then show that these approximate gluing maps can be perturbed to
give a collection of gluing maps $\bga_\al:\sN_\al\to\sC_{W,E}^{*,0}$
and obstruction maps $\bvarphi_\al:\sN_\al\to\sV_\al$ which parametrize
open neighborhoods of the ends of the non-compact moduli space of
$\PU(2)$ monopoles in the following sense: The image $\Imag\bga_\al$ of
a gluing map is a finite-dimensional submanifold of the configuration space
$\sC_{W,E}^{*,0}$ of pairs of connections and spinors; an open
neighborhood $\bga_\al(\bvarphi_\al^{-1}(0))$ in the moduli space
$M_{W,E}^{*,0}$ of $\PU(2)$ monopoles is then cut out of the gluing
map image $\Imag\bga_\al$ by an obstruction section of a finite-rank
obstruction bundle defined over the gluing parameter data $\sN_\al$.

A gluing map $\bga_\al$ is constructed by solving the `infinite-dimensional
part' of the $\PU(2)$ monopole equations \eqref{eq:PT}, essentially
obtained by projecting out the eigenspaces corresponding to the finitely
many `small eigenvalues' tending to zero. More precisely, the scheme we are
forced to use is a variant of that developed by Donaldson
\cite{DonConn,DK}, where we keep the metric fixed and adapt methods of
Taubes \cite{TauFrame,TauStable} to allow us to glue in entire moduli
spaces of anti-self-dual connections on $S^4$: Donaldson's scheme assumes
that the connections are restricted to precompact subsets of their moduli
spaces, while the Riemannian metric on $X$ is allowed to vary
conformally. The obstruction map $\bvarphi_\al$ is then defined by
$\bga_\al$ and the
`finite-dimensional part' of the $\PU(2)$ monopole equations
\eqref{eq:PT} which cannot be solved directly (due to the small
eigenvalues and the resulting growth of Green's operator norms needed to
solve the quasi-linear equation by the Banach space fixed-point
theorem). These small eigenvalues arise here because neither the background
monopole nor the anti-self-dual connections over $S^4$ --- now viewed as
`zero-section $\PU(2)$ monopoles' --- are smooth points of their respective
moduli spaces in the sense of Kodaira-Spencer. These small-eigenvalue
phenomena are reminiscent of those in Taubes' earlier work on gluing
anti-self-dual connections \cite{TauIndef,TauStable} where they arise when
the background connection is trivial. However, for the purposes of
differential-topological calculations, the difficulties surrounding them
can generally be circumvented by working with connections on $\SO(3)$
bundles with non-zero $w_2$ or via blowup tricks \cite{MorganMrowka}: such
a strategy does not work in the case of $\PU(2)$ monopoles.

The construction of gluing and obstruction maps for $\PU(2)$ monopoles
is given in \cite{FL3}, where their existence is established, and the proof
that they parametrize the ends of $M_{W,E}$ is completed in
\cite{FL4}. The difficulties in constructing $\PU(2)$ monopole gluing
maps come from several sources:
\begin{itemize}
\item There are always obstructions to gluing coming from the
anti-self-dual connections over the four-sphere
$S^4$, because of the non-zero cokernel of
the Dirac operator $D_{A}$, and from the background moduli spaces of $\U(1)$
monopoles.
\item The $\PU(2)$ monopole equations, like the Seiberg-Witten
equations, are not conformally invariant. Hence, the
gluing technology for the conformally invariant
anti-self-dual equation developed by Donaldson in
\cite{DonConn,DK} cannot be used directly for $\PU(2)$ monopoles.
\item The gluing theory of Taubes
\cite{TauSelfDual,TauIndef,TauFrame,TauStable} is difficult to adapt to the
case of $\PU(2)$ monopoles because the Bochner
formula for $d_A^+d_A^{+,*}$ --- on which the estimates of
\cite{TauSelfDual,TauIndef,TauFrame,TauStable} rely and which is
well-behaved when the connection $A$ bubbles 
--- must be used in conjunction 
with a Bochner formula for $D_AD_A^*$ which is badly
behaved when the connection $A$ bubbles.
\item In the work of Donaldson \cite{DonConn} and Mrowka
\cite{MrowkaThesis} on the `gluing theorem' for anti-self-dual connections,
the anti-self-dual connections being glued up are assumed to vary
in precompact subsets of their respective moduli spaces. While such
restrictions always simplify the analysis greatly, they cannot be imposed here
since we need to ensure that the entire ends of the moduli
space of $\PU(2)$ monopoles are covered by gluing maps.
\end{itemize}
The Bochner formulas relevant for Taubes' method are given by 
\begin{align*}
2d_A^+d_A^* & = \cov_A^*\cov_A - 2\{\sW^+,\cdot\} + \third R +
\{F_A^+,\cdot\}, \\
D_{A}D_{A}^* &= \cov_A^*\cov_A + \quarter R
+ \half\rho(F(A_{\det W^+}^-)) + \rho(F_A^-).
\end{align*}
The term $F_A^+$ will be uniformly $L^\8$ bounded 
while the term $F_A^-$ is only uniformly bounded in $L^2$ and its $L^\8$ norm
tends to infinity as the connection $A$ bubbles. This phenomenon makes it
extremely difficult to produce Green's operator estimates which are uniform
with respect to a degenerating, approximate $\PU(2)$ monopole $(A,\Phi)$ and
hence solve the equations \eqref{eq:PT} for exact, nearby $\PU(2)$
monopoles. These problems are
overcome in \cite{FL3,FL4} by developing a combination of the gluing
methods of Donaldson and Taubes, but the above difficulties make the gluing
theory and the construction of links much more involved than it is for either 
anti-self-dual connections or Seiberg-Witten monopoles (the simplification in
the latter case stems from the fact that the Seiberg-Witten moduli spaces
are compact \cite{MST}). For example, we need estimates not only for
the gluing maps but also for their differentials (and their inverses) to
prove that the gluing maps are diffeomorphisms and 
cover the moduli space ends \cite{FL4}.

In \cite{FL4} we show that (i) the $\PU(2)$ monopole gluing maps are
`surjective' in the sense that every $\PU(2)$ monopole lies in the image
of a gluing map (so it can be `unglued'), (ii) they are diffeomorphisms
onto their images, and (iii) the gluing map images have an invariant
characterization in the quotient. The surjectivity property of Taubes'
gluing maps for anti-self-dual connections is a special case of a more
general gluing result for critical points of the Yang-Mills functional
\cite[Proposition 8.2]{TauFrame}. Like the proof of a particular case of the
surjectivity statement for anti-self-dual connection gluing maps given by
Donaldson and Kronheimer in \cite[\S 7.2]{DK}, Taubes' argument essentially
relies on estimates for the inverse of the differential of the gluing map
and the `method of continuity' to show that a given point lies in the image
of a gluing map. Again, the main new difficulty here lies in getting
estimates which are uniform with respect to an approximate $\PU(2)$
monopole connection which is `bubbling' (and thus approaching the Uhlenbeck
boundary). Our construction in \cite{FL3,FL4} shows that open neighborhoods
of the lower-level strata of $\barM_{W,E}$ are modelled by zero sets of
sections of finite-rank obstruction bundles: this generalizes the
description given in \S \ref{subsec:LinkOfReducible} of open neighborhoods
of the singular strata in the top level $M_{W,E}$.

\subsection{General reduction formulas and the $\PU(2)$-monopole
analogue of the Kotschick-Morgan conjecture}
\label{subsec:PTConjecture}
In this section we sketch some of the ideas underlying our approach to
the $\PU(2)$-monopole analogue of the Kotschick-Morgan conjecture.

The first observation one needs 
in order to appreciate why the $\PU(2)$-monopole
program should work is that, as discussed in
\S \ref{sec:Cohomology} and shown in \cite{FL2}, the intersection
$\barV(z)\cap \barW(x^{n_a-1})$ of geometric representatives is a
collection of smooth one-manifolds, with one set of boundaries near the
moduli space $M_E^{\asd}$ of anti-self-dual solutions and the other
boundaries in neighborhoods of Seiberg-Witten reducible solutions of the
form 
\begin{equation}
M^{\red}_{W,E_{-\ell},L_1}\times\Sym^\ell(X) \subset IM_{W,E}.
\label{eq:LowerLevelReducibles}
\end{equation}
Because of the obstructions to gluing, it is not clear that all the points of 
\eqref{eq:LowerLevelReducibles} are necessarily contained in
$\barM_{W,E}$, and so $\barM_{W,E}$ may be a proper subset of
$IM_{W,E}$. 

In \cite{FL2} we analyze the intersection of these geometric
representatives in a neighborhood of the anti-self dual solutions and
reducible $\PU(2)$ monopoles in the top Uhlenbeck level (as described here
in \S \ref{sec:Cohomology}).  To generalize Theorem
\ref{thm:CompactReductionFormula} to the case when there are reducible
pairs in the lower levels of $\barM_{W,E}$, we need a precise construction
of the links $\bL_{W,E_{-\ell},L_1}$ of the lower-level reducibles
\eqref{eq:LowerLevelReducibles}.  In \cite{FL4} we use the gluing and
obstruction maps to construct an open neighborhood
$U_{W,E_{-\ell},L_1}$ of the points
\eqref{eq:LowerLevelReducibles} in $\barM_{W,E}$ with a `piecewise
smoothly-stratified boundary'
$$
\bL_{W,E_{-\ell},L_1} := \rd U_{W,E_{-\ell},L_1}.  
$$
This boundary serves as a link of the reducible solutions
\eqref{eq:LowerLevelReducibles} in the
compactified moduli space $\barM_{W,E}$.  Because there are obstructions to
gluing coming from both the background $\PU(2)$ monopoles and the
anti-self-dual connections over $S^4$, it is not known if the Uhlenbeck
compactification has locally finite topology at points in the lower levels.
Although the link given by $\rd U_{W,E_{-\ell},L_1}$ might not have
finite topology, its intersection with the geometric representatives of the
cohomology classes is finite as this intersection takes place in the top
stratum (in the top level, away from any reducibles).

The above remarks suggest that one should obtain a `reduction formula',
conjectured by Pidstrigach and Tyurin, expressing the Donaldson invariants
in terms of integrals over links of Seiberg-Witten moduli spaces:

\begin{conj}[Pidstrigach and Tyurin]\label{conj:ReductionFormula}
If $z\in\AAA(X)$, then
\begin{align*}
2^{n_a-1}D^{c_1(E)}_X(z)
&=
\sum_{L_1\in\overline\Red(W,E)}
\barV(z)\cap\barW(x^{n_a-1})\cap\bL_{W,E_{-\ell},L_1},
\quad\text{if } \deg z = 2d_a, \\
0
&=
\sum_{L_1\in\overline\Red(W,E)}
\barV(z)\cap\barW(x^{n_a-1})\cap\bL_{W,E_{-\ell},L_1},
\quad\text{if } \deg z > 2d_a.
\end{align*}
\end{conj}

Note that the level index $\ell$ appearing in the right-hand side the above
formulas is determined by the reduction $E_{-\ell} = L_1\oplus (\det
E)\otimes L_1$ defined by $L_1$, since $\det E_{-\ell} = \det E$ is fixed
and $c_2(E_{-\ell}) = c_2(E)-\ell$.

The second formula, while not directly interesting, could be useful
in deriving recursion relations determining the intersections
with $\bL_{W,E_{-\ell},L_1}$.
An important step towards proving Witten's conjecture would be to show that
the intersection on the right has some universal expression (whose precise
form might not be known) in terms of Seiberg-Witten invariants:

\begin{conj}[Pidstrigach and Tyurin]\label{conj:PTConjecture}
The pairing on the right-hand side of Conjecture
\ref{conj:ReductionFormula} is given by a universal formula
depending only on $\ell$, $F$, $L_1$, $\SW(\fs_0\otimes L_1)$, 
the intersection form $Q_X$, and invariants of the homotopy type of $X$.
\end{conj}

This is the Pidstrigach-Tyurin version of the `Kotschick-Morgan conjecture'
\cite[Conjecture 6.2.1 \& 6.2.2]{KoM}. More specifically, one would like to
show that the pairing on right-hand side of Conjecture
\ref{conj:ReductionFormula} is given by 
$$
q_X(\ell,F,L_1,Q_X)\cdot\SW(\fs_0\otimes L_1)
$$ 
for some universal polynomial $q_X(\cdot)$, where the
dependence on $X$ is just through its homotopy type (although even getting
the terms on the right-hand side of Conjecture \ref{conj:ReductionFormula}
to be divisible by $\SW(\fs_0\otimes L_1)$ is a highly non-trivial problem).
Naturally, the ultimate aim is to evaluate these pairings explicitly,
following the example of G\"ottsche in
\cite{Goettsche} for the $b^+=1$ wall-crossing formula, and show that they
coincide with the prediction of Witten in the case of simple type.  We gave
calculations of this type for top level reducibles in Theorem
\ref{thm:DegreeZeroFormula}, when $\ell=0$, and
outline the idea for lower-level reducibles below, when $\ell>0$.

The calculations are simplest when $M^{\red}_{W,E_{-\ell},L_1}$ is
zero-dimensional, 
$$
M^{\red}_{W,E_{-\ell},L_1}
= \{[A_r,\Phi_r]\}_{r=1}^n,
$$
so we sketch the basic idea for this special case below. Note that
when $X$ has Seiberg-Witten simple type it may still have
positive-dimensional Seiberg-Witten moduli spaces and though the
associated Seiberg-Witten invariants will vanish, one cannot {\em a
priori\/} rule out their contributions to the Donaldson polynomials. Hence,
even assuming $X$ has Seiberg-Witten simple type,
we still need the thickened moduli spaces of \S
\ref{subsec:LinkOfReducible} to show that positive-dimensional
Seiberg-Witten moduli spaces do not in fact contribute to the Donaldson
polynomials.

Let $\{U_r\}_{r=1}^n$, be neighborhoods of zero in $H^1_{A_r,\Phi_r}$ for
the reducibles $\{[A_r,\Phi_r]\}_{r=1}^n$ in the background moduli space
$M_{W,E_{-\ell}}$ and let $\Gl(U_r,\Sigma)$ be the gluing data associated
with $U_r$ and a (precompact open subset of a)
smooth stratum $\Sigma\subset\Sym^\ell(X)$.  We can cover
a neighborhood of $[A_r,\Phi_r]\times\Sym^\ell(X)$ in $\barM_{W,E}$ with
the images under the gluing maps
$$
\bga_{r,\Sigma}(\bvarphi_{r,\Sigma}^{-1}(0)\cap\Gl(U_r,\Sigma))
$$ 
of the zero loci of the obstruction sections $\bga_{r,\Sigma}$. The pairing
on the right-hand side of Conjecture
\ref{conj:ReductionFormula} then takes the form 
\begin{equation}
\sum_{r=1}^n \barV(z)\cap\barW(x^{n_a-1})\cap\bL_r,
\label{eq:LinkSum}
\end{equation}
where $\bL_r$ is the link of $[A_r,\Phi_r]\times\Sym^\ell(X)$ in
$\cup_\Sigma\bga_{r,\Sigma}(\Gl(U_r,\Sigma))$.
If one could show that the pairing 
$\barV(z)\cap\barW(x^{n_a-1})\cap\bL_r$ were a
multiple of $\sign[A_r,\Phi_r]$, with coefficient independent of
$r$ --- that is, independent of the background pair, then the sum
\eqref{eq:LinkSum} would be a multiple of
$$
\SW(\fs_0\otimes L_1) 
= 
\# M^{\red}_{W,E_{-\ell},L_1} 
=
\sum_{r=1}^n \sign[A_r,\Phi_r].
$$ 
Independence of the background pair can be shown by direct calculation when
$\ell=1$, much as in \cite{Kotschick,KoM,Y}. The fact that the individual
pairings may depend on the background pairs is essentially because the
gluing maps do not quite `commute': gluing up the same gluing data in different
orders yields slightly different composite gluing maps. Similar
difficulties have been encountered in attempts to prove the Kotschick-Morgan
conjecture of Donaldson theory \cite{KoM,MorganOzsvath}.

In the positive dimensional case there are additional problems due to
`spectral flow' or `jumping lines' and this makes it difficult to describe
the links of the lower-level moduli space of $\U(1)$ monopoles,
$M^{\red}_{W,E_{-\ell},L_1}\times\Sym^\ell(X)$. In general, there is no global
Kuranishi model for $M^{\red}_{W,E_{-\ell},L_1}$ which is defined naturally by
small-eigenvalue cutoffs which we can glue up with $S^4$ gluing data
to form open neighborhoods in $M_{W,E}$ --- one encounters
`jumping lines' as the points in a neighborhood of the
background moduli space $M_{W,E_{-\ell}}$ vary. (Models which are global with
respect to the background Seiberg-Witten moduli space are desirable for the
purposes of calculating Euler classes of the obstruction bundles.)  As
outlined in \S
\ref{sec:Cohomology}, we employ stabilization methods \cite{ASFamily,DK} to
address these problems when they are caused by reducibles in the top level
in \cite{FL2}, where no gluing is needed. In the general case, we use
gluing to parametrize the links of lower-level reducibles in combination
with this stabilization procedure \cite{FL3,FL4} 
when $\ell > 0$ and verify Conjectures
\ref{conj:ReductionFormula} and \ref{conj:PTConjecture} 
by direct calculation when $\ell=1$.



\end{document}